\documentclass[12pt]{amsart}

\title[Force to change large cardinal strength]{Force to change large cardinal strength}

\author[Carmody]{Erin Carmody}
\usepackage{amssymb}
\usepackage[hidelinks]{hyperref}
\usepackage{mathrsfs}
\usepackage{setspace}
\usepackage{textcomp}

\usepackage{setspace}
\usepackage[inner=1in, outer=1in, bottom=1in]{geometry}
\include{MathMacrosJDH}

\newtheorem{theorem}{Theorem}
\newtheorem{lemma}[theorem]{Lemma}
\newtheorem{corollary}[theorem]{Corollary}
\newtheorem{definition}[theorem]{Definition}

\begin{document}
\renewcommand{\thepage}{\roman{page}}
\topskip0pt
\vspace*{\fill}
\begin{center}

{ \Huge Force to change large cardinal strength}\\
{ \large by\\
Erin Carmody}\\
\vspace{.5 in}
\begin{singlespace}
{ \large A dissertation submitted to the Graduate Faculty in Mathematics in partial fulfillment of the requirements for the degree of Doctor of Philosophy, The City University of New York.}\\
\vspace{.3 in}
2015
\end{singlespace}
\end{center}
\vspace*{\fill}

\newpage
\renewcommand{\thepage}{\roman{page}}
\topskip0pt
\vspace*{\fill}
{\large
\begin{center}
\textcopyright 2015\\

\vspace{.3 in}

Erin Carmody\\

\vspace{.3 in}

All Rights Reserved
\end{center}
}
\vspace*{\fill}

\newpage
\renewcommand{\thepage}{\roman{page}}
\newpage

\newpage

\newpage
\renewcommand{\thepage}{\roman{page}}
\setcounter{page}{4}
\begin{center}
\textsc{Abstract}
\end{center}

This dissertation includes many theorems which show how to change large cardinal properties with forcing.  I consider in detail the degrees of inaccessible cardinals (an analogue of the classical degrees of Mahlo cardinals) and provide new large cardinal definitions for degrees of inaccessible cardinals extending the hyper-inaccessible hierarchy.  I showed that for every cardinal $\kappa$, and ordinal $\alpha$, if $\kappa$ is $\alpha$-inaccerssible, then there is a $\mathbb{P}$ forcing that $\kappa$ which preserves that $\alpha$-inaccessible but destorys that $\kappa$ is $(\alpha+1)$-inaccessible. I also consider Mahlo cardinals and degrees of Mahlo cardinals. I showed that for every cardinal $\kappa$, and ordinal $\alpha$, there is a notion of forcing $\mathbb{P}$ such that $\kappa$ is still  $\alpha$-Mahlo in the extension, but $\kappa$ is no longer $(\alpha +1)$-Mahlo.  I also show that a cardinal $\kappa$ which is Mahlo in the ground model can have every possible inaccessible degree in the forcing extension, but no longer be Mahlo there.  The thesis includes a collection of results which give forcing notions which change large cardinal strength from weakly compact to weakly measurable, including some earlier work by others that fit this theme. I consider in detail measurable cardinals and Mitchell rank. I show how to change a class of measurable cardinals by forcing to an extension where all measurable cardinals above some fixed ordinal $\alpha$ have Mitchell rank below $\alpha.$  Finally, I consider supercompact cardinals, and a few theorems about strongly compact cardinals.  Here, I show how to change the Mitchell rank for supercompactness for a class of cardinals.

\begin{center}

{ Force to change large cardinal strength}\\
{ \large by\\
Erin Carmody}\\
\end{center}

\vspace{.25 in}

Advisor: Joel David Hamkins\\

This dissertation includes many theorems which show how to change large cardinal properties with forcing.  I consider in detail the degrees of inaccessible cardinals (an analogue of the classical degrees of Mahlo cardinals) and provide new large cardinal definitions for degrees of inaccessible cardinals extending the hyper-inaccessible hierarchy.  I showed that for every cardinal $\kappa$, and ordinal $\alpha$, if $\kappa$ is $\alpha$-inaccerssible, then there is a $\mathbb{P}$ forcing that $\kappa$ which preserves that $\alpha$-inaccessible but destorys that $\kappa$ is $(\alpha+1)$-inaccessible. I also consider Mahlo cardinals and degrees of Mahlo cardinals. I showed that for every cardinal $\kappa$, and ordinal $\alpha$, there is a notion of forcing $\mathbb{P}$ such that $\kappa$ is still  $\alpha$-Mahlo in the extension, but $\kappa$ is no longer $(\alpha +1)$-Mahlo.  I also show that a cardinal $\kappa$ which is Mahlo in the ground model can have every possible inaccessible degree in the forcing extension, but no longer be Mahlo there.  The thesis includes a collection of results which give forcing notions which change large cardinal strength from weakly compact to weakly measurable, including some earlier work by others that fit this theme. I consider in detail measurable cardinals and Mitchell rank. I show how to change a class of measurable cardinals by forcing to an extension where all measurable cardinals above some fixed ordinal $\alpha$ have Mitchell rank below $\alpha.$  Finally, I consider supercompact cardinals, and a few theorems about strongly compact cardinals.  Here, I show how to change the Mitchell rank for supercompactness for a class of cardinals.

\newpage
\renewcommand{\thepage}{\roman{page}}
\begin{center}
\textsc{Aknowledgements}
\end{center}

First of all, I would like to thank Arthur Apter for helping me see the future.  His encouragement has inspired me to see how bright our future is together.  I have total faith and trust in Arthur.  We work well together. It makes me feel good when we can talk about interesting problems together, and I look forward to all of our theorems of the future!\\

From the beginning until the end of this process, Philipp Rothmaler has been a great mentor to me.  He encouraged me to present at seminar, and helped me fall in love with logic. Thank you Philipp for letting me tell you my ideas and for always being nearby to see me through until the end.  .\\

Gunter, thank you for being my friend and my mentor through this process. You have always been there for me, and showed your appreciation of my interests.  Your mind is amazing, and I always learn so much from your talks.  I would like to know more, and create notions of forcing together.\\

Victoria, dear Victoria! You inspire me and you helped me so much.  You took the time to look at every word of my thesis and applications.  You became my sister over time from the oral exam until now where we have our work together. I pray that I can find your kind of focus.  Thank you, Victoria.\\

Dear Roman, thank you for encouraging me in every way possible.  Thank you for being a reliable source of support, and for showing your appreciation. Thank you for everything, Roman.\\

Joel is never-waivering and ever-patient.  He is a great advisor. To see him teach is to watch lightning in a bottle, like an angel defending the line between heaven and hell.  He is enthusiastic and gives 100$\%$ all the time.  His brain moves fast like a computer and his mind is as gentle as a breeze.  He is kind and aware.  He helped me learn how to be a person through his actions and words.  He is so thoughtful and so fun! He is truly-deeply-eternally-vastly, $\dots$ and so on.\\

I would also like to thank my pre-committee: Kaethe Minden, Miha Habic, and Kameryn Williams.  Thank you for forcing me to come to your wonderful seminar to get ready for the final examination.  I could not have done it without your thoughtful comments and friendship all through my time as a student.  You all are brilliant and I know that you are going to create beautiful theorems and change the future of set theory and mathematics.  
\newpage
\renewcommand{\thepage}{\roman{page}}
\tableofcontents

\newpage
\renewcommand{\thepage}{\arabic{page}}
\setcounter{page}{1}
\section{\textbf{Introduction}}

Large cardinals are infinite numbers, so big that we cannot prove their existence in ZFC.  If $\kappa$ is a large cardinal, then $V_{\kappa} \models $ZFC, so proving the existence of a large cardinal would violate Godel's incompleteness theorem.  The theorems in the following chapters assume the consistency of ZFC with large cardinals.  Forcing was developed by Cohen in 1963, and has been used as a way to prove independence results.  Forcing in the following chapters is used to create an extension where one or many large cardinals has some greatest desired degree of some large cardinal property.  Suppose $\kappa \in V$ is a cardinal with large cardinal property $A$.  The purpose of this thesis is to find a notion of forcing $\mathbb{P}$ such that if $G \subseteq \mathbb{P}$ is $V$-generic, the cardinal $\kappa$ no longer has property $A$ in $V[G]$, but has as many large cardinal properties below $A$ as possible.  The main theorems of this dissertation are positive results like this where the large cardinal properties of a cardinal change between the ground model and the forcing extension.  A notion of forcing $\mathbb{P}$ which can change a large cardinal, like the ones in the main theorems, does so carefully.  That is, $\mathbb{P}$ does not just destroy large cardinal properties, it also preserves large cardinal properties.  The project is to pick a desired large cardinal degree and design a notion of forcing which forces a cardinal to lose all large cardinal properties above this degree and keep all large cardinal properties below.  In other words, we are killing large cardinals as softly as possible.  In particular we are killing large cardinals by just one degree from the ground model to the forcing extension.  Thus, a part of this thesis is also devoted to defining degrees of various large cardinal properties such as inaccessible, Mahlo, measurable, and supercompact.

The first chapter is about inaccessible cardinals.  An inaccessible cardinal $\kappa$ is an uncountable cardinal having two properties exhibited by $\omega$, the first infinite cardinal.  Namely, an inaccessible cardinal $\kappa$ is regular: cof$(\kappa) = \kappa$, just like $\omega$ since no sequence of finite length can be cofinal in $\omega$, and an inaccessible cardinal $\kappa$ is also a strong limit: $\forall \beta < \kappa, 2^{\beta} < \kappa$ just like $\omega$ since $\forall n < \omega, 2^n < \omega$.  But inaccessible cardinals are uncountable, unlike $\omega$.  An inaccessible cardinal might also be a limit of inaccessible cardinals, this is a 1-inaccessible cardinal.  A cardinal $\kappa$ is $\alpha$-inaccessible if and only if $\kappa$ is inaccessible and for every $\beta < \alpha$, $\kappa$ is a limit of $\beta$-inaccessible cardinals. The main theorem of the first chapter is that if $\kappa$ is $\alpha$-inaccessible, there is a forcing extension where $\kappa$ is $\alpha$-inaccessible but not $(\alpha + 1)$-inaccessible.  If a cardinal $\kappa$ is $\kappa$-inaccessible it is called hyper-inaccessible.  In the first chapter, I define $\alpha$-hyper$^{\beta}$-inaccessible cardinals and beyond this to the richly-inaccessible cardinals, utterly-inaccessible cardinals, deeply-inaccessible cardinals, truly-inaccessible cardinals, eternally-inaccessible cardinals, and vastly-inaccessible cardinals.  And, I develop a notation system for a hierarchy of inaccessibility reaching beyond this.  

The second chapter is about forcing to change degrees of Mahlo cardinals.  A Mahlo cardinal $\kappa$ has many inaccessible cardinals below, every closed unbounded (club) subset of $\kappa$ contains an inaccessible cardinal.  Let $\mathcal{I}(X)$ be an operator which gives the set of inaccessible limit points of a set $X$.  A cardinal $\kappa$ is called greatly inaccessible if and only if there is a uniform, normal filter on $\kappa$, closed under $\mathcal{I}(X)$.  The first theorem of the second chapter shows that a greatly inaccessible cardinal is equivalent to a Mahlo cardinal.  Since greatly inaccessible cardinals are every possible inaccessible degree, as defined in chapter 1, Mahlo cardinals are every possible inaccessible degree defined.  One of the main theorems of the second chapter is to force a Mahlo cardinal $\kappa$ to lose its Mahlo property in the extension while preserving that $\kappa$ is every inaccessible degree defined.  There is also a classical hierarchy of Mahlo cardinal degrees, like the hierarchy of inaccessible cardinal degrees defined in the first chapter.  A cardinal $\kappa$ is 1-Mahlo if and only if the set of Mahlo cardinals below $\kappa$ is stationary in $\kappa$ (every club in $\kappa$ contains a Mahlo cardinal).  A cardinal $\kappa$ is $\alpha$-Mahlo if and only if $\kappa$ is Mahlo and for all $\beta < \alpha$ the set of $\beta$-Mahlo cardinals is stationary in $\kappa$.  The main theorem of the second chapter is that for any $\kappa$ and $\alpha$ where $\kappa$ is $\alpha$-Mahlo in $V$ there is a forcing extension where $\kappa$ is $\alpha$-Mahlo but not $(\alpha +1)$-Mahlo.

The third chapter involves many large cardinal notions.  A few of the theorems in chapter three are about weakly compact cardinals.  A cardinal $\kappa$ is weakly compact if and only if $\kappa$ is inaccessible and for every tree of height $\kappa$ whose levels have cardinality less than $\kappa$ has a branch of length $\kappa$.  One theorem shows that weakly compact cardinals are every classical Mahlo degree.  Another shows that one can force a weakly compact cardinal $\kappa$ to lose its weak compactness in the extension while keeping all of its Mahlo degrees.  This chapter includes theorems by others, especially Victoria Gitman about weakly compact cardinals and the following large cardinal properties which are defined in chapter three:  weakly measurable cardinals, strongly Ramsey cardinals, and ineffable cardinals.

The fourth chapter is about measurable cardinals, which have all of the large cardinal properties described so far.  A cardinal $\kappa$ is measurable if and only if it is the critical point of an elementary embedding $j:V \to M$.  The classic hierarchy of degrees for measurable cardinals is the Mitchell rank, denoted $o(\kappa)$ for a cardinal $\kappa$, is defined in chapter four.  In this chapter we are able to force to change the Mitchell rank for a class of cardinals in the ground model.  The main theorem is that for any $V \models$ ZFC, any $\alpha \in$ Ord, and any $\kappa > \alpha$ in $V$, there is a forcing extension $V[G]$ where $o(\kappa) = \min\{\alpha, o(\kappa)^V\}$.  In other words, in the forcing extension every large cardinal above $\alpha$ has Mitchell rank at most $\alpha$.

The fifth and final chapter of this dissertation is about supercompact and strongly compact cardinals.  A cardinal $\kappa$ is $\theta$-supercompact if and only if there is an elementary embedding $j:V \to M$ with critical point $\kappa$ and $M^{\theta} \subseteq M$.  That is, the target of the elementary embedding is closed under $\theta$-sequences.  Equivalently, $\kappa$ is $\theta$-supercompact if and only if there is a normal fine measure on $P_{\kappa}\theta$.  A cardinal $\kappa$ is $\theta$-strongly compact if and only if there is a $\kappa$-complete fine measure on $P_{\kappa} \theta$.  The first main results of the fifth chapter are due to a key result by Hamkins and Shelah and show how to force a supercompact (or strongly compact) cardinal to be at most $\theta$-supercompact (or $\theta$-strongly compact) for a fixed $\theta$. Also, Magidor's result that one can force to separate strongly compact and supercompact cardinals is in this chapter.  In chapter five, Mitchell rank for $\theta$-supercompact cardinals, denoted $o_{\theta \text{-sc}}(\kappa)$, for a cardinal $\kappa$ is discussed in detail.  The main result shows how to find a forcing extension where for many $\kappa, \alpha, \theta$  in $V$ the Mitchell rank for $\theta$-supercompactness for $\kappa$ is at most $\alpha$ in the extension: $o_{\theta \text{-sc}}(\kappa)^{V[G]} = \min\{ \alpha, o_{\theta \text{-sc}}(\kappa)^V \}$.

\section{\textbf{Degrees of Inaccessible Cardinals}}

In this section, I have found forcing extensions where the degree of an inaccessible cardinal is reduced to an exact specified amount from its degree in the ground model.  In general terms, for a given cardinal $\kappa$, we find a notion of forcing which simultaneously preserves large cardinal properties of $\kappa$ up to some level while destroying its large cardinal properties above that level. Suppose $\kappa$ is a cardinal with large cardinal property $A$. Thus, $\kappa$ also has property $B$ for any lower property $B$ which follows from $A$.  The objective is to find a forcing notion $\mathbb{P}$, with generic object $G \subseteq \mathbb{P}$, such that in the forcing extension, $V[G]$, the cardinal $\kappa$ no longer has property $A$, but still has property $B$.  The main theorem below shows how to do this when the difference between $A$ and $B$ is one degree of inaccessibility.   An infinite cardinal $\kappa$ is \textit{inaccessible} if and only if $\kappa$ is an uncountable, regular, strong limit cardinal.  A cardinal $\kappa$ is \textit{1-inaccessible} if and only if $\kappa$ is inaccessible and a limit of inaccessible cardinals.  A cardinal $\kappa$ is $\alpha$-\textit{inaccessible} if and only if $\kappa$ is inaccessible, and for every $\beta < \alpha$, the cardinal $\kappa$, is a limit of $\beta$-inaccessible cardinals. 

\begin{theorem} \label{Theorem. Change one inaccessible degree.}

If $\kappa$ is $\alpha$-inaccessible, then there is a forcing extension where $\kappa$ is still $\alpha$-inaccessible, but not $(\alpha + 1)$-inaccessible. 
\end{theorem}

The following Lemma establishes some basic facts about degrees of inaccessible cardinals. 

\begin{lemma} \label{Lemma. Basic facts of inaccessible degrees.}

1. A cardinal $\kappa$ is 0-inaccessible if and only if $\kappa$ is inaccessible.  \\
2. If $\kappa$ is $\alpha$-inaccessible, and $\beta < \alpha$, then $\kappa$ is also $\beta$-inaccessible.\\
3. A cardinal $\kappa$ cannot be $\eta$-inaccessible, for any $\eta > \kappa$.
\begin{proof}
1.  If $\kappa$ is 0-inaccessible, then by definition $\kappa$ is inaccessible and a limit of $\beta$-inaccessible cardinals for every $\beta < 0$, hence $\kappa$ is inaccessible.  If $\kappa$ is inaccessible then it is vacuously true that for every $\beta < 0$, the cardinal $\kappa$ is a limit of $\beta$-inaccessible cardinals.  

2. Suppose $\kappa$ is $\alpha$-inaccessible, and $\beta < \alpha$. Since $\kappa$ is $\alpha$-inaccessible, for every $\eta < \alpha$, the cardinal $\kappa$ is a limit of $\eta$-inaccessible cardinals.  Since every $\gamma$ less than $\beta$ is also less than $\alpha$, the cardinal $\kappa$ is a limit of $\gamma$-inaccessible cardinals, for every $\gamma < \beta$.  Thus, $\kappa$ is also $\beta$-inaccessible.

3. By way of contradiction, suppose $\kappa$ is the least cardinal with the property that $\kappa$ is $(\kappa + 1)$-inaccessible. It follows from the definition that $\kappa$ is a limit of $\kappa$-inaccessible cardinals. Thus, there is a cardinal $\beta < \kappa$ for which $\beta$ is $\kappa$-inaccessible. Since $\kappa$ is a limit ordinal, $\beta + 1 < \kappa$, and so by statement (2), $\beta$ is also $(\beta + 1)$-inaccessible.  This contradicts that $\kappa$ is the least cardinal with this property. It follows from the lemma that $\kappa$ cannot be $\eta$-inaccessible for any $\eta \ge \kappa + 1$, because being $\eta$-inaccessible would imply that $\kappa$ is $(\kappa + 1)$-inaccessible.
\end{proof}
\end{lemma}

And now the proof of Theorem \ref{Theorem. Change one inaccessible degree.}:

\begin{proof}

Let $\kappa$ be $\alpha$-inaccessible. If $\kappa$ is not $(\alpha + 1)$-inaccessible, then trivial forcing will give the forcing extension where $\kappa$ is $\alpha$-inaccessible, but not $(\alpha + 1)$-inaccessible. Thus, assume that $\kappa$ is $(\alpha + 1)$-inaccessible, and we will find a forcing extension where it is no longer $(\alpha + 1)$-inaccessible but still $\alpha$-inaccessible. The idea of the proof is to add a club, $C$, to $\kappa$ which contains no $\alpha$-inaccessible cardinals, and then force to change the continuum function to kill strong limits which are not limit points of $C$.  To change the continuum function, we will perform Easton forcing. 

Easton forcing to change the continuum function works when the GCH holds in the ground model. For our purposes, we only need the GCH pattern to hold up through $\kappa$, by forcing with $\mathbb{P}$, which is a $\kappa$-length iteration of $Add(\gamma,1)$ for regular $\gamma \in V^{\mathbb{P}_{\gamma}}$.  The forcing $\mathbb{P}$, neither destroys, nor creates, inaccessible cardinals below $\kappa$, and thus is a mild preparatory forcing, so that we are in a place from which we will eventually perform Easton forcing.  

First, see that $\mathbb{P}$ does not destroy inaccessible cardinals, by investigating its factors. Let $\eta < \kappa$ be inaccessible, and let $p \in \mathbb{P}$. Let $p^{\le \beta}$ be the condition, $p$, restricted to domain $[0, \beta]$.  And, let $p^{> \beta}$, be the condition, $p$, with domain restricted to $(\beta, \kappa]$.  Let $\mathbb{P}^{\le \beta} = \{ p^{\le \beta} \ : \ p \in \mathbb{P} \}$ and $\mathbb{P}^{> \beta} = \{ p^{< \beta} \ : \ p \in \mathbb{P} \}$, so that $\mathbb{P}$ factors as $\mathbb{P}^{\le \beta} \ast \mathbb{P}^{ > \beta}$.  For $\beta < \eta$, the first factor, $\mathbb{P}^{\le \beta}$, is small relative to $\eta$, and thus cannot force $2^{\beta} \ge \eta$.  And, the second factor, $\mathbb{P}^{ > \beta}$, is ${\le}\beta$-closed, so adds no new subsets to $\beta$. Thus $\eta$ is still a strong limit after forcing with $\mathbb{P}$.  Similar factoring arguments show that $\eta$ remains a regular limit after forcing with $\mathbb{P}$ [Jech, 232+].  Since $\mathbb{P}$ preserves inaccessible cardinals, it follows by induction, that it preserves any degree of inaccessibility defined so far.  Further, since inaccessibility is downward absolute, it follows by induction, that $\mathbb{P}$ cannot create $\beta$-inaccessible cardinals for any $\beta < \kappa$.  Let $G \subseteq \mathbb{P}$ be $V$-generic, and force to $V[G]$, where $V[G] \models GCH$, and where $V$ and $V[G]$ have the same $\beta$-inaccessible cardinals for any $\beta < \kappa$.  

Next, force with $\mathbb{C}$, which will add a club subset to $\kappa$, which contains no $\alpha$-inaccessible cardinals in the forcing extension.  Conditions $c \in \mathbb{C}$ are closed, bounded, subsets of $\kappa$, consisting of infinite cardinals and containing no $\alpha$-inaccessible cardinals.  The forcing $\mathbb{C}$ is ordered by end-extension: $d \le c$ if and only if $c = d \cap (\sup(c) + 1)$.  Let $H \subseteq \mathbb{C}$ be $V[G]$-generic, and let $C = \cup H$.  Then, it is claimed that $C$ is club in $\kappa$, and contains no $\alpha$-inaccessible cardinals in $V[G][H]$.  First, let us see that $C$ is unbounded.  Let $\beta \in \kappa$, and let $D_{\beta} = \{ d \in \mathbb{C} \ : \ \max(d) > \beta \}$.  The set $D_{\beta}$ is dense in $\mathbb{C}$, since, if $c \in \mathbb{P}$, the set $d = c \cup \{ \beta' \}$, where $\beta' > \beta$ and $\beta'$ is not $\alpha$-inaccessible, is a condition in $D_{\beta}$ and $d$ is stronger than $c$.  Hence, there is $c_{\beta} \in D_{\beta} \cap G$, so that there is a an element of $C$ above $\beta$.  Thus, in $V[G][H]$, the set $C$ is unbounded in $\kappa$.  Suppose $\delta \cap C$ is unbounded in $\delta < \kappa$.  Since $C$ is unbounded, there is $\delta' \in C$, where $\delta' > \delta$, and a condition $c_{\delta} \in C$ which contains $\delta'$.  Since the conditions are ordered by end-extension, and since there is a sequence (possibly of length 1) of conditions which witness that $\delta \cap C$ is unbounded in $\delta$, which all must be contained in the condition $c_{\delta}$ (which contains an element above $\delta$), it follows that $\delta \cap c_{\delta}$ is unbounded in $\delta$.  Since $c_{\delta}$ is closed, $\delta \in c_{\delta}$.  Therefore, $\delta \in C$, which shows that $C$ is also closed.   

The forcing, $\mathbb{C}$, adds a new club, $C$, to $\kappa$.  It remains to show that $C$ contains no $\alpha$-inaccessible cardinals, $C$ contains unboundedly many $\beta$-inaccessible cardinals, and that $\mathbb{C}$ preserves cardinals, cofinalities and strong limits.  The new club, $C$, does not contain any ground model $\alpha$-inaccessible cardinals, since if it did, there would be a condition in $G$ which contains an $\alpha$-inaccessible, contradicting that no condition in $G$ contains an $\alpha$-inaccessible.  To see that the new club, $C$, contains unboundedly many ground model $\beta$-inaccessible cardinals, for every $\beta < \alpha$, fix $\beta < \alpha$ and $\eta < \kappa$.  Let $D_{\eta}$ be the set of conditions in $\mathbb{C}$ which contain a $\beta$-inaccessible $\gamma$ above $\eta$, and which contain a sequence of inaccessible cardinals, unbounded in $\eta$, witnessing that $\gamma$ is $\beta$-inaccessible.  Let us see that $D_{\eta}$ is dense in $\mathbb{C}$.  Let $c \in \mathbb{C}$.  Let $\gamma$ be the next $\beta$-inaccessible above both $\eta$ and and the maximal element of $c$.  Since $\gamma$ is the next $\beta$-inaccessible past $\eta$ and $\beta < \alpha$, there are no $\alpha$-inaccessible cardinals in $(\eta, \gamma]$.  Also, this block, $(\eta, \gamma]$, contains the tails of all sequences of inaccessible cardinals, which witnesses that $\gamma$ is $\beta$-inaccessible.  Let $d = c \cup \{ (\eta, \gamma] \} \cap $CARD. Then, $d$ extends $c$, contains a $\beta$-inaccessible above $\eta$, and a sequence of inaccessible cardinals which witness that $\gamma$ is $\beta$-inaccessible.  Thus, $d \in D_{\eta}$, which shows that $D_{\eta}$ is dense in $\mathbb{C}$, and thus shows that $C$ contains unboundedly many $\beta$-inaccessible cardinals.  From this fact also follows, that, in the final extension, $\kappa$ is still $\alpha$-inaccessible, as we shall soon see.  Finally, the forcing $\mathbb{C}$ preserves cardinals and cofinalities greater than or equal to $\kappa + 1$, forcing over $V[G] \models$ GCH, since $|\mathbb{C}| = \kappa^{ < \kappa } = \kappa$.  

For $\beta < \kappa$, the forcing, $\mathbb{C}$, is not ${\le}\beta$-closed, since if $\beta$ is $\alpha$-inaccessible, there is a $\beta$-sequence of conditions unbounded in $\beta$, but no condition could close the sequence since it would have to include $\beta$.  However, for every $\beta < \kappa$, the set $D_{\beta} = \{ d \in \mathbb{C} \ : \ \max(d) \ge \beta \}$ is dense in $\mathbb{C}$ and is ${\le}\beta$-closed.  This is true, since for any $\beta$-sequence of conditions in $D_{\beta}$, one can close the sequence by taking unions at limits and adding the top point, which cannot be inaccessible, because it is not regular, since this top point is above $\beta$, but has cofinality $\beta$.  Thus, for every $\beta < \kappa$, the forcing, $\mathbb{C}$, is forcing equivalent to $D_{\beta}$, which is $\le \beta$-closed.  Thus, $\mathbb{C}$ preserves all cardinals, cofinalities, and strong limits.  Thus. $V[G][H] \models$ GCH and we have forced to add $C \subseteq \kappa$, club, which contains no $\alpha$-inaccessible cardinals. 

The last step of the proof is to force over $V[G][H]$, with $\mathbb{E}$, Easton's forcing to change the continuum function.  Specifically, let $\mathbb{E}$ force $2^{\gamma^+} = \delta^+$, where $\gamma \in C$, for $\gamma$ infinite, and $\delta$ is the next element of $C$ past $\gamma^+$.  This forcing preserves all cardinals and cofinalities [Jech, 232+], and also preserves that $\kappa$ is inaccessible since factoring at any $\beta < \kappa$ reveals that the first factor is too small to force $2^{\beta}$ up to $\kappa$, and the second factor adds no new subsets to $\beta$, thus preserving that $\kappa$ is a strong limit cardinal, in addition to preserving that $\kappa$ is a regular limit cardinal.  However, $\mathbb{E}$ does not preserve all inaccessible cardinals below $\kappa$.  In fact, $\mathbb{E}$ destroys all strong limits which are not limit points of $C$.  Let $\eta < \kappa$ be a strong limit which is not in $C'$, the set of limit points of $C$.  Since $C$ is closed, there is a greatest element of $C$ below $\eta$, call it $\gamma$.  Then, the next element of $C$, call it $\delta$, is greater than or equal to $\eta$.  So, when $\mathbb{E}$ forces $2^{\gamma^+} = \delta^+$, it destroys that $\eta$ is a strong limit since $\gamma < \eta$ implies $\gamma^+ < \eta$ (since $\eta$ was a strong limit in $V[G][H]$).  Let $K \subseteq \mathbb{E}$ be $V[G][H]$-generic.  Since $C'$ contains no ground model $\alpha$-inaccessible cardinals, there are no $\alpha$-inaccessible cardinals below $\kappa$ in $V[G][H][K]$.  Thus, $\kappa$ is not $(\alpha + 1)$-inaccessible in the final forcing extension.  

Finally, see that since $C'$ contains unboundedly many ground model $\beta$-inaccessible cardinals (which is the same as the set of $\beta$-inaccessible cardinals in the intermediate extensions since $\mathbb{P}$ and $\mathbb{C}$ preserve all inaccessible cardinals) for every $\beta < \alpha$.  The cardinal $\kappa$ is still $\beta$-inaccessible in $V[G][H][K]$.  Let $\gamma' < \gamma$, where $\gamma \in C'$ and $\gamma$ is inaccessible.  Since $\gamma$ is a limit point of $C$, the next element of $C$, call it $\delta$, which is above $\gamma'$, is also below $\gamma$, and $\delta^+ < \gamma$ since $\gamma$ is a strong limit.  Thus, $2^{\gamma'} \le 2^{\gamma'^+} = \delta^+ < \gamma$.  So, $\gamma$ is still a strong limit.  Since $\mathbb{E}$ preserves cardinals and cofinalities, this shows that $\gamma$ is still inaccessible in the final extension.  Thus, all inaccessible limit points of $C$ are preserved.  The earlier density argument also shows that for any $\beta < \alpha$ it is dense that $C$ contains an interval containing all cardinals and also contains a $\beta$-inaccessible cardinal.  Since all limit points of such an interval are limit points of $C$, it follows that all inaccessible cardinals in such an interval are preserved.  Thus, all degrees of inaccessible cardinals are preserved. If not, then suppose the least degree of inaccessibility which is destroyed an $\eta$-inaccessible cardinal for some $\eta < \alpha$.  Then for some $\gamma < \eta$, there is no sequence of $\gamma$-inaccessible cardinals witnessing this $\eta$-inaccessible cardinal, which implies that $\gamma$-inaccessible cardinals are destroyed, which contradicts that $\eta$ was the least degree not preserved.  Thus, in $V[G][H][K]$, for every $\beta < \alpha$, there are unboundedly many $\beta$-inaccessible cardinals below $\kappa$.  Thus, $\kappa$ is still $\alpha$-inaccessible in the final extension. But, $\kappa$ is not $(\alpha + 1)$-inaccessible in $V[G][H][K]$ since $C'$ contains no $\alpha$-inaccessible cardinals.

\end{proof}

Theorem \ref{Theorem. Change one inaccessible degree.} shows how to change the degree of an inaccessible cardinal, in a forcing extension, to be exactly $\alpha$ when $\kappa$ has degree at least $\alpha$ in the ground model. The following theorem shows how to force to a universe where there are no inaccessible cardinals, but where every every ground model weakly inaccessible cardinal is still weakly inaccessible.

\begin{theorem} \label{Theorem. Distinguish class of inaccessible from weakly inaccessible.}

For any $V \models ZFC$ there exists $V[G]$ with no inaccessible cardinals, but where every ground model weakly inaccessible cardinal is still weakly inaccessible. 

\begin{proof}
The proof is very similar to the proof of Theorem \ref{Theorem. Change one inaccessible degree.}.  It is a class forcing, where the first step is to force that Ord is not Mahlo.  We add a club $C \subseteq Ord$ such that $C$ contains no inaccessible cardinals.  The second step is to perform Easton forcing to make $2^{\gamma^+} = \delta^+$ whenever $\gamma \in C$ and where $\delta$ is the next element of $C$ above $\gamma$. The idea of the proof appears as a footnote in Hamkins[5].  The proof of Theorem \ref{Theorem. Change one inaccessible degree.} shows that the combined forcing preserves all strong limits of $C$ which are internally strong limits of $C$, and destroys all strong limits of $C$ which are not internally strong limits of $C$.  Since $C$ contains no inaccessible cardinals, there are no inaccessible cardinals in the final extension.  However, as one can see from the proof of Theorem \ref{Theorem. Change one inaccessible degree.} or from Lemma \ref{Lemma. Club shooting.} that the forcings preserve all cardinals and cofinalities, so that all regular limit cardinals are preserved.  Hence all weakly compact inaccessible cardinals are preserved.  
\end{proof}
\end{theorem}

 If a cardinal $\kappa$ is $\kappa$-inaccessible, then it is defined to be \textit{hyper-inaccessible}. Therefore, if $\kappa$ is hyper-inaccessible, Theorem \ref{Theorem. Change one inaccessible degree.} shows how to force to make $\kappa$ have inaccessible degree $\alpha$ for some $\alpha < \kappa$.  One can force to change a hyper-inaccessible cardinal to have maximal degree by forcing to add a club which avoids all degrees above the desired degree, and then force with Easton forcing to destroy strong limits which are not limit points of the new club, as in the proof of Theorem \ref{Theorem. Change one inaccessible degree.}. I'll state it as a corollary below.

\begin{corollary}

If $\kappa$ is hyper-inaccessible, then for any $\alpha < \kappa$, there exists a forcing extension where $\kappa$ is $\alpha$-inaccessible, but not hyper-inaccessible.  

\end{corollary}

Lemma \ref{Lemma. Basic facts of inaccessible degrees.}  shows that the greatest degree of $\alpha$-inaccessibility that $\kappa$ can be is $\kappa$-inaccessible.  Mahlo began the investigation of degrees of inaccessible cardinals [Kan], and fully defined the analgous notions for Mahlo cardinals, and I shall continue his work by formalizing the degrees of inaccessible cardinals in the remainder of this section.  So, how do we proceed beyond hyper-inaccessible cardinals in defining degrees of inaccessibility?  By repeating the process:  a cardinal $\kappa$ is \textit{1-hyper-inaccessible} if and only if $\kappa$ is hyper-inaccessible, and a limit of hyper-inaccessible cardinals (0-hyper-inaccessible is hyper-inaccessible).  That is, $\kappa$ is 1-hyper-inaccessible if and only if the set $\{ \gamma < \kappa \ : \ \gamma \mbox{ is } \gamma\mbox{-inaccessible} \}$ is unbounded in $\kappa$.  A cardinal $\kappa$ is $\alpha$-\textit{hyper-inaccessible} if and only if $\kappa$ is hyper-inaccessible, and for every $\beta < \alpha$, the cardinal $\kappa$ is a limit of $\beta$-hyper-inaccessible cardinals.  A cardinal $\kappa$ is \textit{hyper-hyper-inaccessible}, denoted hyper$^2$-inaccessible, if and only if $\kappa$ is $\kappa$-hyper-inaccessible (hyper$^0$-inaccessible denotes inaccessible).  The following theorem shows we have again reached an apparent roadblock in the hierarchy of inaccessible cardinals.

\begin{theorem}

If $\kappa$ is $\alpha$-hyper-inaccessible, then $\alpha \le \kappa$.  

\begin{proof}

Suppose $\kappa$ is the least cardinal with the property that $\kappa$ is $(\kappa + 1)$-hyper-inaccessible.  Then $\kappa$ is a limit of $\kappa$-hyper-inaccessible cardinals.  Let $\beta < \kappa$ be some $\kappa$-hyper-inaccessible cardinal. Then, for every $\gamma < \kappa$, the cardinal $\beta$ is a limit of $\gamma$-hyper inaccessible cardinals.  This implies that for every $\gamma < \beta + 1 < \kappa$, the cardinal $\beta$ is a limit of $\gamma$-hyper-inaccessible cardinals.  Thus, by definition, the cardinal $\beta$ is $(\beta + 1)$-hyper-inaccessible, which contradicts that $\kappa$ is the least cardinal with this property.

\end{proof}

\end{theorem}

In order to reach higher degrees of inaccessibility past this limit, repeat the process: a cardinal $\kappa$ is \textit{1-hyper$^2$-inaccessible} if and only if $\kappa$ inaccessible and a limit of hyper$^2$-inaccessible cardinals.  A cardinal $\kappa$ is \textit{$\alpha$-hyper$^2$-inaccessible} if and only if, for every $\beta < \alpha$, the cardinal $\kappa$ is inaccessible and a limit of $\beta$-hyper$^2$-inaccessible cardinals.  By the same argument from the previous proof, $\alpha$ must be less than of equal to $\kappa$ in this definition.  And, similarly, it is defined that $\kappa$ is \textit{hyper$^3$-inaccessible} if and only if $\kappa$ is $\kappa$-hyper$^2$-inaccessible.  Thus, we arrive at the beginning of the formulation of the continuation of Mahlo's work with a general definition of degrees of hyper-inaccessibility:

\begin{definition} \label{Definition. Hyper inaccessible cardinals.}

A cardinal $\kappa$ is \textit{$\alpha$-hyper$^{\beta}$-inaccessible} if and only if \\
1) the cardinal $\kappa$ is inaccessible, and\\
2) for all $\eta < \beta$, the cardinal $\kappa$ is $\kappa$-hyper$^{\eta}$-inaccessible, and \\
3) for all $\gamma < \alpha$, the cardinal $\kappa$ is a limit of $\gamma$-hyper$^{\beta}$-inaccessible cardinals.

\end{definition}

Definition \ref{Definition. Hyper inaccessible cardinals.} subsumes the previous definitions of hyper-inaccessible cardinals. For example, 0-hyper$^0$-inaccessible is just inaccessible, since the second and third parts of the definition are not applicable.  If $\kappa$ is $\alpha$-hyper$^0$-inaccessible, then by definition, $\kappa$ is inaccessible and for all $\gamma < \alpha$, the cardinal $\kappa$ is a limit of $\gamma$-hyper$^0$-inaccessible cardinals, hence $\kappa$ is $\alpha$-inaccessible. And if $\kappa$ is $\alpha$-inaccessible, then it is inaccessible and a limit of $\beta$-inaccessible cardinals, for every $\beta < \alpha$, hence $\kappa$ is $\alpha$-hyper$^0$-inaccessible.  And, 0-hyper-inaccessible is just hyper-inaccessible since, from the general definition, this means that $\kappa$ is inaccessible and for every $\eta < 1$, the cardinal $\kappa$ is $\kappa$-hyper$^{\eta}$-inaccessible, i.e., $\kappa$ is $\kappa$-hyper$^{0}$-inaccessible, hence hyper-inaccessible.  In the previous definition of hyper-inaccessible, it was required that $\kappa$ be hyper-inaccessible, but this requirement is included in the second part of Definition \ref{Definition. Hyper inaccessible cardinals.} since this implies that $\kappa$ is hyper-inaccessible whenever $\beta > 0$.

Definition \ref{Definition. Hyper inaccessible cardinals.} gives a general definition for hyper-inaccessible cardinals.  The following theorem shows that a cardinal $\kappa$ can be at most hyper$^{\kappa}$-inaccessible, using this definition.

\begin{theorem}

A cardinal $\kappa$ cannot be 1-hyper$^{\kappa}$-inaccessible.

\begin{proof}

Suppose $\kappa$ is the least cardinal with the property that $\kappa$ is 1-hyper$^{\kappa}$-inaccessible.  Then, by definition, $\kappa$ is a limit of hyper$^{\kappa}$-inaccessible cardinals.  Let $\beta < \kappa$ be hyper$^{\kappa}$-inaccessible.  Then, for all $\gamma < \kappa$, the cardinal $\beta$ is $\beta$-hyper$^{\gamma}$-inaccessible.  In particular, since $\beta < \kappa$, the cardinal $\beta$ is $\beta$-hyper$^{\beta}$-inaccessible.  This implies that $\beta$ is 1-hyper$^{\beta}$-inaccessible, contradicting that $\kappa$ is the least cardinal with this property.  

\end{proof}
\end{theorem}

In order to define degrees of inaccessible cardinals, as in Definition \ref{Definition. Hyper inaccessible cardinals.}, beyond hyper-inaccessible, we need more words.  Here, I will introduce new adjectives to describe higher degrees of inaccessible cardinals, analogous  to hyper-inaccessible degrees.  If $\kappa$ is a hyper$^{\kappa}$-inaccessible cardinal, call it \textit{richly-inaccessible}. Then, proceed as before:  $\kappa$ is \textit{1-richly-inaccessible} if and only if $\kappa$ is inaccessible and a limit of richly-inaccessible cardinals. Then, follow the same pattern as before to keep finding higher degrees past supposed obstacles: if $\kappa$ is $\kappa$-richly-inaccessible, then $\kappa$ is called \textit{hyper-richly-inaccessible}.  If $\kappa$ is $\kappa$-hyper-richly-inaccessible, then $\kappa$ is called \textit{hyper$^2$-richly-inaccessible}.  If $\kappa$ is hyper$^{\kappa}$-richly-inaccessible, then $\kappa$ is \textit{richly$^2$-inaccessible}.  Then, with these words, the greatest $\kappa$ can be is \textit{richly$^{\kappa}$-inaccessible}.  Thus, define $\kappa$ to be \textit{utterly-inaccessible} if and only if $\kappa$ is richly$^{\kappa}$-inaccessible, so that the next degree of inaccessibility is 1-utterly-inaccessible. A cardinal $\kappa$ is richly-utterly-inaccessible if and only if $\kappa$ is hyper$^{\kappa}$-utterly-inaccessible. A cardinal $\kappa$ is utterly$^2$-inaccessible if and only if $\kappa$ is richly$^{\kappa}$-utterly-inaccessible.  This goes on forever:  if a cardinal $\kappa$ is utterly$^{\kappa}$-inaccessible it is called \textit{deeply-inaccessible}; if a cardinal $\kappa$ is deeply$^{\kappa}$-inaccessible, it is called \textit{truly-inaccessible}; if a cardinal $\kappa$ is truly$^{\kappa}$-inaccessible, it is called \textit{eternally-inaccessible}; if a cardinal $\kappa$ is eternally$^{\kappa}$-inaccessible, it is called \textit{vastly-inaccessible}, and so on.  

I will end this section with an exploration of how to characterize these degrees.  Let $I$ be the class of inaccessible cardinals.  Let $\mathcal{I}$ be the \textit{inaccessible limit point operator}, a class operator which assigns to a class, its inaccessible limit points:

$$ \mathcal{I}(X) = \{ \alpha \ : \ \mbox{ $\alpha$ is an inaccessible limit point of $X$} \} $$

And, for classes $X_{\alpha}$, let $\Delta_{\alpha \in ORD} X_{\alpha} = \{ \gamma \ : \ \gamma \in \cap_{\alpha < \gamma} X_{\alpha} \}$ denote the diagonal intersection of a collection (of classes).  All of the degrees of inaccessible cardinals defined so far are also definable from $I, \mathcal{I},$ and $\Delta$ as follows. First, define $\mathcal{I}^{\alpha} (X)$ for all $\alpha$:

$$ \mathcal{I}^0(X) = X $$
$$ \mathcal{I}^{\alpha + 1}(X) = \mathcal{I}(\mathcal{I}^{\alpha}(X))$$
$$ \mathcal{I}^{\lambda}(X) = \bigcap_{\beta < \lambda} \mathcal{I}^{\beta}(X)$$

Then, $\mathcal{I}^0(I) = I$, the class of inaccessible cardinals.  And, $\mathcal{I}(I) = \mathcal{I}^1(I)$ is the class of 1-inaccessible cardinals, $\mathcal{I}^2(I)$ is the class of 2-inaccessible cardinals, and so forth so that $\mathcal{I}^{\alpha}(I)$ is the class of $\alpha$-inaccessible cardinals. Thus, the class of hyper-inaccessible cardinals is 

$$ \Delta_{\alpha \in \text{Ord}}(\mathcal{I}^{\alpha}(I)) = \{ \kappa \ : \ \kappa \in \bigcap_{\alpha < \kappa} \mathcal{I}^{\alpha}(I) \}$$

Let $H = \Delta_{\alpha \in ORD} (\mathcal{I}^{\alpha}(I))$ be the class of hyper-inaccessible cardinals.  Let $H_{1} = \mathcal{I}(H)$ be the class of 1-hyper-inaccessible cardinals.  The class of $\alpha$-hyper-inaccessible cardinals is denoted $H_{\alpha}$. If we take the diagonal intersection of these classes we get $H^2 = \Delta_{\alpha \in ORD}(H_{\alpha})$, the class of hyper$^2$-inaccessible cardinals. Then, applying the inaccessible limit point operator to this class gives $\mathcal{I}(H^2) = H_1^2$, the class of 1-hyper$^2$-inaccessible cardinals.  Then $H^3 = \Delta_{\alpha \in ORD} (H_{\alpha}^2)$, the class of hyper$^3$-inaccessible cardinals. Then, taking the diagonal intersection of all the $H^{\alpha}$ degrees of hyper-inaccessibility gives $R = \Delta_{\alpha \in ORD} (H^{\alpha})$, the class of richly-inaccessible cardinals.  Apply $\mathcal{I}$ to $R$ to get $R_1, R_2, \dots, R_{\alpha}$, the class of $\alpha$-richly-inaccessible cardinals. Applying diagonal intersection to this class gives $HR = \Delta_{\alpha \in ORD} (R_{\alpha})$, the class of hyper-richly-inaccessible cardinals. Apply $\mathcal{I}$ to this class to get $HR_1 = \mathcal{I}(HR)$, the class of 1-hyper-richly-inaccessible cardinals. Then $H^2R = \Delta_{\alpha \in ORD}(HR_{\alpha})$, the class of hyper$^2$-richly-inaccessible cardinals.  Proceed in this way to get $H^{\alpha}R$, the class of hyper$^{\alpha}$-richly-inaccessible cardinals. Apply the diagonal intersection to these classes to get $R^2 = \Delta_{\alpha \in ORD} (H^{\alpha}R)$, the class of richly$^2$ inaccessible cardinals.  Build up the classes, $R^{\alpha}$, of richly$^{\alpha}$-inaccessible cardinals and take the diagonal intersection of these to get $U = \Delta_{\alpha \in ORD} (R^{\alpha})$, the class of utterly-inaccessible cardinals. Then $HRUDTEV$ is the class of hyper-richly-utterly-deeply-truly-eternally-vastly-inaccessible cardinals, and so on.

The words hyper, richly, utterly, and so on mark places in the process of defining degrees of inaccessible cardinals where we take a diagonal intersection over $\alpha \in$ Ord of a collection of classes.  Thus, we have designed a notation system for the classes of inaccessible cardinals using meta-ordinals.  This notation system captures all the classes described so far, and can go as far as you like.  It is a notation system for meta-ordinals which is like Cantor's normal form for ordinals, but instead of $\omega$ we use $\Omega$, a symbol for the order-type of Ord.  The degree of any inaccessible cardinal will be denoted by $t$, a formal syntactic expression for a meta-ordinal of the form $\Omega^{\alpha}\cdot \beta + \Omega^{\eta} \cdot \gamma + \dots + \Omega \cdot \delta + \sigma$ where $\alpha > \eta > \dots  \in $Ord, and $\beta, \gamma, \delta, \sigma \in$ Ord.  If $s$ and $t$ are meta-ordinals then the ordering is essentially lexicographical.  If the degree on $\Omega$ of the leading term of $t$ is greater than $s$, then $s < t$.  If $s$ and $t$ have the same greatest degree of $\Omega$, then compare the coefficients of the leading term.  If these are the same, compare the next highest powers of $\Omega$ in $s$ and $t$ and so on.  Now that we have these meta-ordinals, we can describe the classes of degrees of inaccessible cardinals in a uniform way.  If $\eta$ and $\kappa$ are hyper-inaccessible cardinals then $\eta$ is in the class of $\eta$-inaccessible cardinals and $\kappa$ is in the class of $\kappa$-inaccessible cardinals.  But now I can describe both as being in the class of $\Omega$-inaccessible cardinals, which is the diagonal intersection of all classes of $\alpha$-inaccessible cardinals, for $\alpha \in$ Ord.  Defined in this way, the degree of an inaccessible cardinal $\kappa$ can be described as $t$-inaccessible for some meta-ordinal $t$.  The only restriction being that all of the ordinals in $t$ are less than or equal to $\kappa$.  In this way, inaccessible cardinals with the same degree of inaccessibility can be described with the same meta-ordinal.  

\begin{definition} \label{Definition. t-inaccessible.}
If $t$ is a meta-ordinal term with ordinals less than $\kappa$, then a cardinal $\kappa$ is $t$-inaccessible if and only if $\kappa$ is inaccessible and for every meta-ordinal term $s < t$ having ordinals less than $\kappa$, the cardinal $\kappa$ is a limit of $s$-inaccessible cardinals. 
\end{definition}

This definition of $t$-inaccessibility for a cardinal $\kappa$ can be extended to include meta-ordinal terms $t$ with ordinals less than or equal to $\kappa$ by replacing $\kappa$ by $\Omega$.  For example, a cardinal $\kappa$ which is $\kappa$-inaccessible is now called $\Omega$-inaccessible. A cardinal $\kappa$ which is $\Omega \cdot \kappa$-inaccessible is a hyper$^{\kappa}$-inaccessible cardinal which is better called a richly-inaccessible cardinal, or an $\Omega^2$-inaccessible cardinal.

 Here are a few of the classes described previously with this notation system:

\noindent $\kappa$ is $\Omega$-inaccessible $\iff$ $\kappa$ is hyper-inaccessible\\
$\kappa$ is $\Omega^2$-inaccessible $\iff$ $\kappa$ is richly-inaccessible\\
$\kappa$ is $(\Omega^2 + \Omega)$-inaccessible $\iff$ $\kappa$ is hyper-richly-inaccessible\\
$\kappa$ is $\Omega^3$-inaccessible $\iff$ $\kappa$ is utterly-inaccessible\\
$\kappa$ is $(\Omega^3 + \Omega^2)$-inaccessible $\iff$ $\kappa$ is richly-utterly-inaccessible\\
$\kappa$ is $\Omega^3 \cdot 2$-inaccessible $\iff$ $\kappa$ is utterly$^{2}$-inaccessible\\
$\kappa$ is $\Omega^4$-inaccessible $\iff$ $\kappa$ is deeply-inaccessible\\
$\kappa$ is $\Omega^5$-inaccessible $\iff$ $\kappa$ is truly-inaccessible\\
$\kappa$ is $\Omega^6$-inaccessible $\iff$ $\kappa$ is eternally-inaccessible\\
$\kappa$ is $\Omega^7$-inaccessible $\iff$ $\kappa$ is vastly-inaccessible\\
$\kappa$ is $(\Omega^7 + \Omega^6 + \Omega^5 + \Omega^4 + \Omega^3 + \Omega^2  + \Omega + \alpha)$-inaccessible $\iff$ $\kappa$ is $\alpha$-hyper-richly-utterly-deeply-truly-eternally-vastly-inaccessible\\

Now that we have this uniform notation, we have the following theorem:

\begin{theorem} If $\kappa$ is $t$-inaccessible, then there is a forcing extension where $\kappa$ is still $t$-inaccessible, but not $(t+1)$-inaccessible.

\begin{proof} Sketch.  The proof is very similar to the proof of Theorem  \ref{Theorem. Change one inaccessible degree.}.  Assume $V \models$ ZFC. Suppose $\kappa \in V$ is $t$-inaccessible, where $t$ is a meta-ordinal term with only ordinal exponents and coefficients less than $\kappa$. Force with $\mathbb{C}$. 
Conditions are closed, bounded subsets of $\kappa$ such that if $c \in \mathbb{C}$ then $c$ contains 
no $t$-inaccessible cardinals. The forcing $\mathbb{C}$ preserves cardinals and cofinalities.
Since $\forall \gamma < \kappa$, the set $D_{\gamma} = \{ d \in \mathbb{C} : \max(d) > \gamma \}$ is dense, the new club $C = \cup G \subseteq \mathbb{C}$ is unbounded.
It follows that $C$ is closed, since it is unbounded and the conditions are closed. For every $ s < t$, where $s$ is a meta-ordinal term with ordinal parameters less than $\kappa$, the set $D_s = \{ d \in \mathbb{C} \ : \ d \text {contains a block of cardinals up to an $s$-inaccessible cardinal} \}$ is dense:
if $c \in \mathbb{C}$, then if $\eta$ is the least $s$-inaccessible cardinal above $\eta$, then
$d = c \cup ( \{ [\max(C), \eta] \} \cap \text{CARD} )$ is in $D_s$.\\

Next force with $\mathbb{E}$: For infinite $\gamma \in C$, force $2^{\gamma^+} = \delta^+$ where $\delta$ is the next element of $C$ past $\gamma^+$.
$\mathbb{E}$ destroys strong limits which are not in $C'$:
If $\eta$ is a strong limit which is not a limit point of $C$, then there is a maximal element of $C$ below $\eta$.  Call it $\gamma$.  
Then $\gamma^+ < \eta$ and the next element of $C$, call it $\delta$, is above $\eta$.  Thus $2^{\gamma^{+}} = \delta^+ > \eta$ destroys
that $\eta$ is a strong limit.  If $\eta \in C'$ and $\eta$ is a strong limit, then $\delta^+ < \eta$ so $2^{\gamma^{+}} < \eta$ for every $\gamma < \eta$.
Thus, $\eta$ is still a strong limit in the final extension.\\

Since $C$ contains unboundedly many ground model $s$-inaccessible cardinals for $s < t$, with a block of cardinals below them, then in the final extension these cardinals are still $s$-inaccessible (by a proof by meta-ordinal induction on $s$), and so $\kappa$ is still $t$-inaccessible. However since all $t$-inaccessible cardinals below $\kappa$ are no longer strong limits, $\kappa$ is not $t+1$-inaccessible.
\end{proof}
\end{theorem}

Finally, consider the connection between sets of inaccessible cardinals below $\kappa$ and regressive functions on $\kappa$.  Consider the set $\mathcal{I}^2(I)$ of 2-inaccessible cardinals. If $\kappa \notin \mathcal{I}^2(I)$ this means that there is a tail $[\beta, \kappa)$ where, for each $\gamma$ inaccessible in $[\beta, \kappa)$, there is, $b_{\gamma}$, the greatest inaccessible cardinal below $\gamma$.  Thus, I can define a regressive function $f: [\beta, \kappa) \to \kappa$ with $f(\alpha) = \beta_{\alpha}$.  Thus, for every $\alpha$, there are no inaccessible cardinals between $f(\alpha)$ and $\alpha$. In general, I can say that $\kappa \notin \mathcal{I}^2(I)$ if and only if there is a regressive function whose domain is a tail of $\kappa$, with the property that there are no inaccessible cardinals between $f(\alpha)$ and $\alpha$ for any $\alpha$ in the domain.  Similarly, I can define $\kappa \notin \mathcal{I}^3(I)$ if and only if there is a regressive function $f: \kappa \to \kappa$ such that for all inaccessible $\alpha < \kappa$, we have $\mathcal{I}^2(I) \cap (f(\alpha), \alpha) = \emptyset$, and there exists a regressive function $g_{\alpha} : \beta \to \beta$ such that for all $\beta \in (f(\alpha), \alpha)$, we have $\mathcal{I}^1(I) \cap (g_{\alpha}(\beta), \beta) = \emptyset$.  Thus, no $\alpha < \kappa$ is 2-inaccessible, and thus $\kappa$ cannot be a limit of 2-inaccessible cardinals so $\kappa \notin \mathcal{I}^3(I)$, but $\kappa$ could still be 2-inaccessible.

\newpage

\section{\textbf{Mahlo Cardinals}}

This section begins with theorems about Mahlo cardinals, showing that Mahlo cardinals have all the inaccessible large cardinal properties from the last section.  An analogue of the classical notion of greatly Mahlo is defined for inaccessible cardinals. Also, the classical degrees of Mahlo cardinals are described, and theorems follow which show how to force a large cardinal $\kappa$ to have a desired Mahlo degree in the extension.   An infinite cardinal $\kappa$ is \textit{Mahlo} if and only if the set of inaccessible cardinals below $\kappa$ is a stationary subset of $\kappa$.  A cardinal $\kappa$ is \textit{greatly inaccessible} if and only if there is a uniform, normal filter on $\kappa$, closed under the inaccessible limit point operator:

$$\mathcal{I}(X) = \{ \alpha \in X \ : \ \alpha \mbox{ is an inaccessible limit point of } X \}$$. 

A uniform filter on $\kappa$ has for every $\beta < \kappa$, the set $[\beta, \kappa)$ is in the filter, and normal means that the filter is closed under diagonal intersections $\Delta_{\alpha < \kappa}$. The last part of the definition means that if $X$ is in the filter, so is $\mathcal{I}(X)$.  The first theorem of this section shows that greatly inaccessible is equivalent to Mahlo.  

\begin{theorem} \label{Theorem. Mahlo same as greatly inaccessible.}
A cardinal $\kappa$ is greatly inaccessible if and only if $\kappa$ is Mahlo.   

\begin{proof}
For the forward implication, suppose $\kappa$ is Mahlo.  Let $F$ be the filter generated by sets of the form $C \cap I$, where $C$ is club in $\kappa$, and $I$ is the set of inaccessible cardinals below $\kappa$.  Note that $\kappa$ is club in itself, and $\kappa \cap I = I$, so $I \in F$.  Then, the claim is that $F$ is a uniform, normal filter, closed under $\mathcal{I}$.  First, if $C, D$ are club in $\kappa$, then $C \cap D$ is club. Thus, for any clubs $C$ and $D$, the equation $(C \cap I) \cap (D \cap I) = (C \cap D) \cap I$ implies that any set $A$ in the filter, generated by sets of the form $C \cap I$, where $C$ is club, is itself a superset of a set of the form $E \cap I$, where $E$ is club in $\kappa$.  Then, $\emptyset \notin F$, since the empty set has no non-trivial subsets, hence cannot be a superset of the form $C \cap I$.  Next, if $A \in F$, and $A \subseteq B$, then there is a club $C \in F$ such that $C \cap I \subseteq A \subseteq B$, thus $B \in F$, by construction, since $F$ is the filter generated by sets of this form.  Third, if $A$ and $B$ are elements of $F$, then there are clubs $C$ and $D$ such that $C \cap I \subseteq A$, and $D \cap I \subseteq B$, so $A \cap B$ contains $(C \cap D) \cap I$. This is of the form which generated the filter, thus $A \cap B \in F$.  Also, $F$ is uniform since if $b$ is a bounded subset of $\kappa$, then $\kappa \setminus b$ contains a club $E$.  Thus, $E \cap I \subseteq E \subseteq \kappa \setminus b$, so $\kappa \setminus b \in F$.  Since the cardinality of any co-bounded set is $\kappa$, the filter is uniform.  It remains to show that $F$ is normal, and closed under the inaccessible limit point operation.  To see that $F$ is normal, i.e. closed under diagonal intersection, see first that club sets are closed under diagonal intersection.  Let $C_{\alpha}$ be club for every $\alpha < \kappa$.  Then $\Delta C_{\alpha}$ is closed, since, if $\delta$ is a limit point of $\Delta C_{\alpha}$, then $\delta$ is a limit point of each $C_{\alpha}$, for $\alpha < \delta$.  To see this, let $\alpha < \delta$.  Then, since $\delta$ is a limit point of $\Delta C_{\alpha}$, there exists a $\gamma$ such that $\alpha < \gamma < \delta$, with $\gamma \in \Delta C_{\alpha}$, i.e., $\gamma \in \cap_{\beta < \gamma} C_{\beta}$. Since $\alpha < \gamma$, this implies that $\gamma \in C_{\alpha}$.  Hence, that $\delta$ is a limit point of $\Delta C_{\alpha}$ implies $\delta$ is a limit point of $C_{\alpha}$, for all $\alpha < \delta$.  Thus, since each $C_{\alpha}$ is closed, $\delta \in C_{\alpha}$, for all $\alpha < \delta$.  Thus $\delta \in \cap_{\alpha < \delta} C_{\alpha}$.  Thus, $\delta \in \Delta C_{\alpha}$.  Thus, $\Delta C_{\alpha}$ is closed.  To see that $\Delta C_{\alpha}$ is unbounded, observe that for any $\lambda < \kappa$, the set $\cap_{\alpha < \lambda} C_{\alpha}$ is club.  I'll show now that $\Delta C_{\alpha}$ is unbounded, using this fact.  Suppose that $\beta < \kappa$.  Since $\cap_{\alpha < \beta} C_{\alpha}$ is unbounded, there is $\beta_1 \in \cap_{\alpha < \beta} C_{\alpha}$, with $\beta_1 > \beta$.  Then since $\cap_{\alpha < \beta_1} C_{\alpha}$ is unbounded, there is $\beta_2 \in \cap_{\alpha < \beta_1} C_{\alpha}$ with $\beta_2 > \beta_1$.  Continue in this way, to get an increasing sequence $\langle \beta_{n} \rangle$, where $\beta_{n+1} \in \cap_{\alpha < \beta_n} C_{\alpha}$. Taking the limit of this sequence, $\lim_{n \to \infty} \beta_n = \beta_{\omega}$, and the limit of the intersection of clubs gives that $\beta_{\omega} \in \cap_{\alpha < \beta_{\omega}} C_{\alpha}$, and $\beta_{\omega} > \beta$, since it is the limit of an increasing sequence above $\beta$.  Thus, $\beta \in \Delta C_{\alpha}$, is a witness that $\Delta C_{\alpha}$ is unbounded.  Finally, for this direction, to show Mahlo implies greatly inaccessible, $F$ is closed under $\mathcal{I}$, since if $A \in F$, then $\mathcal{I}(A) = I \cap A' \cap A$, and $A'$ is club, so $I \cap A' \in F$, and $A \in F$, so that $I \cap A' \cap A \in F$, and thus $\mathcal{I}(A) \in F$. 

For the other direction, if $\kappa$ is greatly inaccessible, then the uniform, normal filter, $F$, on $\kappa$, contains all the club subsets of $\kappa$.  This is true since all tails $(\beta, \kappa)$ are in $F$, since $\kappa$ is regular, and $F$ is uniform, and all clubs can be written as the diagonal intersection of tails.  Let $C$ be club in $\kappa$.  Then define

$$
A_{\alpha}^C =
\begin{cases}
[\alpha, \kappa), & \text{if }\alpha \in C \\
[\beta, \kappa), & \text{if }\alpha \notin C, \mbox{ where } \beta \text{ next element of }C, \mbox{ past } \alpha
\end{cases}
$$

It is claimed that $C = \Delta A_{\alpha}^C$.  Let $x \in C$.  Then, for every $\gamma < x$, the set $A_{\gamma}^C$ contains $[x, \kappa)$, since the next element of $C$ past $\gamma$ is at most $x$.  Thus, $x \in \cap_{\gamma < x} A_{\gamma}^C$ so that $x \in \Delta_{\gamma < x} A_{\gamma}^C$.  For the other direction, suppose $x \in \Delta A_{\alpha}^C$.  Then, for every $\gamma < x$, $x \in A_{\gamma}^C$. Thus, let $\gamma < x$, then either $\gamma \in C$ or the next element of $C$ past $\gamma$ is below $x$ since $x \in A_{\gamma}^C$.  Thus, below $x$, there is an unbounded subset of $x$ whose elements are in $C$.  Thus, $x$ is a limit point of $C$.  Since $C$ is closed, $x \in C$.  Thus $C = \Delta A_{\alpha}^C$.  Thus, since $F$ is normal, $C \in F$.  Then, since $\kappa \in F$, and $\mathcal{I}(\kappa) = I \in F$, and for any club $D \subseteq \kappa$, the set $D$ is in $F$, and their intersection, $D \cap I$, is in $F$, hence $D \cap I$ is non-empty.  Thus $I$ is stationary in $\kappa$, and thus $\kappa$ is Mahlo.  

\end{proof}

\end{theorem}

Theorem \ref{Theorem. Mahlo every inaccessible degree.} below shows that a Mahlo cardinal is every degree of inaccessibility defined previously.  However, a cardinal being every degree of inaccessibility is not equivalent to the cardinal being Mahlo.  Theorem \ref{Theorem. Distinguish Mahlo and all inaccessible degrees.} separates the two notions with forcing to destroy the Mahlo property of a cardinal, while preserving that the cardinal is every inaccessible degree.  

\begin{theorem} \label{Theorem. Mahlo every inaccessible degree.}
If $\kappa$ is Mahlo, then for every meta-ordinal term $t$ having ordinals less than $\kappa$, the cardinal $\kappa$ is a $t$-inaccessible cardinal as in Definition \ref{Definition. t-inaccessible.})

\begin{proof}
If $\kappa$ is Mahlo, then $\kappa$ is greatly inaccessible, by Theorem \ref{Theorem. Mahlo same as greatly inaccessible.}.  Hence, there is a uniform, normal filter $F$ on $\kappa$, closed under the inaccessible limit point operation. Since $\kappa \in F$, and $\mathcal{I}(\kappa) = I \in F$, all sets of $\alpha$-inaccessible cardinals for $\alpha < \kappa$, below $\kappa$, are in $F$.  And, since $F$ is normal, all sets of hyper inaccessible degrees below $\kappa$, which are all definable from diagonal intersection, are in $F$.  Hence, every set of inaccessible cardinals definable from $I, \mathcal{I},$ and $\Delta$, in the manner of the previous section, by closing under these operations, is unbounded in $\kappa$.  
\end{proof}

\end{theorem}
Theorem \ref{Theorem. Change one inaccessible degree.} of the last section shows that any cardinal which is at least $\alpha$-inaccessible in the ground model, for some $\alpha$, can be made to be no more than $\alpha$-inaccessible in a forcing extension.  Therefore, if $\kappa$ is $\alpha$-inaccessible, then there is a forcing extension where $\kappa$ is $\alpha$-inaccessible, but not Mahlo.  The following theorem shows that if $\kappa$ is Mahlo, there is a forcing extension where $\kappa$ is no longer Mahlo but still any degree of inaccessibility defined in the previous section.

\begin{theorem} \label{Theorem. Distinguish Mahlo and all inaccessible degrees.}

If $\kappa$ is Mahlo, then there is a forcing extension where $\kappa$ is $t$-inaccessible for every meta-ordinal term $t$ with ordinals less than or $\kappa$ , but where $\kappa$ is not Mahlo.

\begin{proof}  
Let $\kappa$ be a Mahlo cardinal.  Let $\mathbb{C}$ be the forcing to add a club, $\mathcal{C} = \cup G$, where $G \subseteq \mathbb{C}$ is $V$-generic, which contains no inaccessible cardinals.  The argument in Theorem \ref{Theorem. Change one inaccessible degree.} shows that $\mathbb{C}$ preserves cardinals, cofinalities, and all inaccessible cardinals in $V_{\kappa}$.  In fact, $\mathbb{C}$ did not change $V_{\kappa}$ at all, $V_{\kappa}^{V[G]} = V_{\kappa}$.  Since $V_{\kappa}$ has all the sets needed to define that $\kappa$ is $\alpha$-inaccessible, for cardinals $\gamma < \kappa$, being $\alpha$-inaccessible is absolute to $V_{\kappa}$.  Hence $\mathbb{C}$ preserves all inaccessible degrees.  Hence $\kappa$ is still every possible inaccessible degree.  However, $\kappa$ is no longer Mahlo in $V[G]$, since $C \cap I = \emptyset$.  

\end{proof}

\end{theorem}

Just as there are infinitely many degrees of inaccessible cardinals, there are infinitely many Mahlo degrees.  One might be tempted to define the next degree of Mahlo to be a Mahlo limit of Mahlo cardinals, exactly as was done with the degrees of inaccessible cardinals, and indeed there is a hierarchy of Mahlo cardinals that can be defined this way.  However, there is a more powerful and appropriate way to define the degrees of Mahlo cardinals which is much stronger.  The definition of Mahlo cardinal is primarily about stationary sets, and the degrees of Mahlo cardinals are classically defined by stationary sets. Thus, the classical degrees of Mahlo cardinals are defined using stationary sets.  Namely, an infinite cardinal $\kappa$ is \textit{1-Mahlo} if and only if $\kappa$ is Mahlo, and the set of Mahlo cardinals below $\kappa$ is stationary in $\kappa$.  In general, $\kappa$ is \textit{$\alpha$-Mahlo} if and only if $\kappa$ is Mahlo, and for every $\beta < \alpha$, the set of $\beta$-Mahlo cardinals below $\kappa$ is stationary in $\kappa$.  The degrees of Mahlo cardinals go on forever, just as the degrees of inaccessible cardinals; if $\kappa$ is $\kappa$-Mahlo, then $\kappa$ is \textit{hyper-Mahlo}, and so on. 

The main theorem of this section will show how to force to change degrees of Mahlo cardinals by adding a club avoiding a stationary set of cardinals of a certain Mahlo degree, while preserving all stationary subsets of cardinals of a lesser Mahlo degree.  A modification of the forcing $\mathbb{C}$, to add a club, from the proof of Theorem \ref{Theorem. Change one inaccessible degree.}, will work.  Given two sets $A$ and $B$, it is said that $A$ \textit{reflects} $B$ if for some $\delta \in B$, the set $A \cap \delta$ is stationary in $\delta$.  The following lemma shows that if $A$ does not reflect in $B$, then the forcing to add a club through $B$ will preserve the stationary subsets of $B$.  

That is, if we force to add a club which avoids the strong set, then the stationarity of the weak set and all of its stationary subsets will be preserved if the strong set does not reflect in the weak set.  In the proof of the theorem, the strong set will be Mahlo cardinals of a fixed degree, and the weak set will be the set of Mahlo cardinals of a lesser degree which are not of the fixed degree.  

\begin{lemma} \label{Lemma. Weak and strong.}

If $\kappa$ is Mahlo, and the sets $A$ and $B$ partition the inaccessible cardinals below $\kappa$, where $A$ does not reflect in $B$, then the forcing to add a club avoiding $A$ will preserve all the stationary subsets of $B$.  Furthermore, the forcing does not add sets to $V_{\kappa}$.

\begin{proof}

Suppose $S \in V$ is a stationary subset of $B$.  Let $\mathbb{P}$ be the forcing to add a club through the complement of $A$, and let $\dot{E}$ be a name for a club subset of $\kappa$, and $p$ a condition which forces that $E$ is club in $\kappa$. Let $\delta \in S$ be such that $\langle V_{\delta}, \in, A \cap \delta, B \cap \delta, \dot{E} \cap V_{\delta} \rangle$ is an elementary substructure of $\langle V_{\kappa}, \in, A, B, \dot{E} \rangle$.  We can find such $\delta \in S$, since the set of $\delta$ giving rise to elementary substructures is club and $S$ is stationary.  Since $\delta \in B$, by the assumption, there is a club set $c \subseteq \delta$, containing no points from $A$.  Now we construct a pseudo-generic $\delta$-sequence of conditions below $\delta$, deciding more and more about $\dot{E}$, using the elements of $c$ to guide the construction.  Build a descending sequence of conditions $\langle c_{\alpha} \ : \ \alpha < \delta \rangle$ in $\mathbb{P} \cap V_{\delta}$, below $p$, and given $c_{\alpha}$, choose $c_{\alpha + 1}$ to force a specific ordinal above $\alpha$ into $\dot{E}$ with $\sup(c_{\alpha + 1}) \in c$, and at limits $\lambda$, let $c_{\lambda} = (\cup_{\alpha < \lambda} c_{\alpha}) \cup \{ \sup(\cup_{\alpha < \lambda} c_{\alpha} \}$.  Notice that $\sup(\cup_{\alpha < \lambda} c_{\alpha}) \in c$, since $c$ is closed, and therefore is not in $A$, so that $c_{\lambda}$ is a condition for limit ordinals $\lambda$, in the construction.  That is, we can get through the limit steps below $\delta$, precisely because $c$ contains no points from $A$.  Let $c^* = (\cup_{\alpha < \delta}) \cup \{ \delta \}$.  Then $c^* \subseteq \delta$ is club in $\delta$ that decides $\dot{E}$ in a way that is unbounded in $\delta$.  Thus, $c^*$ is a condition which forces $\dot{E}$ meets $S$.  Thus, $S$ must still be stationary in the extension. Finally, this club shooting forcing does not add sets to $V_{\kappa}$ by Lemma \ref{Lemma. Club shooting.}

\end{proof}

\end{lemma}

The main result of this section is that if $\kappa$ is $\alpha$-Mahlo, then there is a forcing extension where $\kappa$ is still $\alpha$-Mahlo, but no longer $(\alpha + 1)$-Mahlo. That is, $\kappa$ is forced to maximal Mahlo degree, $\alpha$-Mahlo.  In fact, a very similar argument can be used to show that if $\kappa$ is Mahlo of any fixed hyper-degree (for example if $\kappa$ is hyper-richly-Mahlo), it can be forced to be $\alpha$-Mahlo (but not hyper-richly Mahlo).

\begin{theorem} \label{Theorem. Change one Mahlo degree.}

If $\kappa$ is $\alpha$-Mahlo, then there is a forcing extension where $\kappa$ is $\alpha$-Mahlo, but not $(\alpha + 1)$-Mahlo.
\begin{proof}

Let $\alpha < \kappa$ be fixed.  Suppose $\kappa$ is $\alpha$-Mahlo.  Let $A$ be the set of $\alpha$-Mahlo cardinals below $\kappa$.  Note that $A$ does not reflect in its complement in the inaccessible cardinals, for if, for some $\beta < \kappa$, inaccessible, in the complement of $A$, the set of $\alpha$-Mahlo cardinals below $\beta$ is stationary, then $\beta$ is $(\alpha + 1)$-Mahlo, hence $\alpha$-Mahlo, hence $\beta \in A$.  But, $\beta$ was supposed to be an inaccessible cardinal in the complement of $A$.  Thus, $A$ does not reflect in the set of inaccessible cardinals in its complement.  

If $\kappa$ is not $(\alpha + 1)$-Mahlo, then force trivially to show the result. Let $\mathbb{C}$ be the notion of forcing which adds a club $\mathcal{C}$ through the complement of $A$.  Conditions are ordered by end-extension.  Then by Lemma \ref{Lemma. Club shooting.} the forcing $\mathbb{C}$ preserves truth in $V_{\kappa}$, preserves cardinals, cofinalities, and adds a club to $\kappa$, disjoint from $A.$  Since $A$ is no longer stationary in $\kappa$, it is no longer $(\alpha + 1)$-Mahlo.  Since $A$ contains no clubs, $\kappa \setminus A$ is stationary.  Thus, by Lemma \ref{Lemma. Weak and strong.}, since $A$ does not reflect in its complement, all stationary subsets of $\kappa \setminus A$ are preserved.  Thus, I only need to show that $\forall \beta < \alpha$ the set $T_{\beta} = \{ \ \delta < \kappa \ : \ \delta \text{ is } \beta \text{-Mahlo but not } (\beta + 1) \text{-Mahlo} \ \}$ is stationary. Let $C \subseteq \kappa$ be club.  Let $\delta$ be least in $C'$ which is $\beta$-Mahlo.  Then $C \cap \delta$ is club in $\delta$ so that $(C \cap \delta)' \cap \delta = C' \cap \delta$ is club in $\delta$ and contains no $\beta$-Mahlo cardinals by the minimality of $\delta$.  Thus $\delta$ is not $(\beta + 1)$-Mahlo.  Thus $T_{\beta}$ is stationary.  Thus the stationarity of $T_{\beta}$ is preserved.  Thus $\kappa$ is still $\alpha$-Mahlo.  

\end{proof}

\end{theorem}

As in the previous section for inaccessible cardinals (Definition \ref{Definition. Hyper inaccessible cardinals.}), we can define the higher degrees of Mahlo cardinals:
\begin{definition} \label{Definition. Hyper Mahlo cardinals.}

A cardinal $\kappa$ is \textit{$\alpha$-hyper$^{\beta}$-Mahlo} if and only if \\
1) the cardinal $\kappa$ is Mahlo, and\\
2) for all $\eta < \beta$, the cardinal $\kappa$ is $\kappa$-hyper$^{\eta}$-Mahlo, and \\
3) for all $\gamma < \alpha$, the set of $\gamma$-hyper$^{\beta}$-Mahlo cardinals below $\kappa$ is stationary in $\kappa$.

\end{definition}

Also, as in Definition \ref{Definition. t-inaccessible.} for degrees of inaccessible cardinals, we can define the classes of degrees of Mahlo cardinals uniformly by using meta-ordinals:

\begin{definition} \label{Definition. t-Mahlo.}
If $t$ is a meta-ordinal term having only ordinals less than $\kappa$, then a cardinal $\kappa$ is $t$-Mahlo if and only if for every meta-ordinal term $s < t$ having only ordinals less than $\kappa$, the set of $s$-inaccessible cardinals is stationary in $\kappa$.  

\end{definition}

\begin{theorem}
If $\kappa$ is $t$-Mahlo, where $t$ is a meta-ordinal term having parameters less than $\kappa$, then there is a forcing extension $V[G]$ where $\kappa$ is $t$-Mahlo, but not $(t+1)$-Mahlo.

\begin{proof}
Sketch.  The proof is a generalization of the proof of Theorem \ref{Theorem. Change one Mahlo degree.}.  The proof is to add a club $C \subseteq \kappa$ which avoids the $t$-Mahlo cardinals below $\kappa$.
\end{proof}
\end{theorem}

To end this section, we take a look at a large cardinal property below inaccessible, and large cardinal properties just above inaccessible.  A cardinal $\kappa$ is \textit{worldly} if and only if $V_{\kappa}$ is a model of ZFC.  If $\kappa$ is inaccessible, then it is worldly since $V_{\kappa} \models$ ZFC.  Therefore, the first theorem is about singular worldly cardinals, which are not inaccessible.  A cardinal $\kappa$ is $\Sigma_{n}$-\textit{worldly} if and only if $V_{\kappa}$ satisfies the $\Sigma_n$ fragment of ZFC.  The following theorem shows how to force a worldly cardinal to be $\Sigma_n$-worldly, but not $\Sigma_{n+1}$-worldly [Hamkins3].

\begin{theorem} (Hamkins) If $\theta$ is any singular worldly cardinal, then for any $n \in \mathbb{N}$, there is a forcing extension where $\theta$ is not worldly, but still $\Sigma_n$-worldly.  
\end{theorem}

The next theorem shows that it is not always possible to force between large cardinal properties.  It shows that if $\kappa$ is inaccessible, there may not be a forcing extension where $\kappa$ is worldly but not inaccessible.  

\begin{theorem}  (Fuchs)  If the existence of an inaccessible cardinals is consistent with ZFC, then it is consistent with ZFC that a cardinal $\kappa$ is inaccessible and has no forcing extension where $\kappa$ is worldly and not inaccessible. 
\end{theorem}
 Let me say a little about the proof.  This proof will show that if $\kappa$ is a regular worldly cardinal, it may not be possible to force that $\kappa$ is a singular worldly cardinal.  Suppose $\kappa$ is a regular worldly cardinal and there is a forcing extension $V[G]$ where $\kappa$ is a singular worldly cardinal.  Then, by a covering lemma argument with the Dodd-Jensen core model, it follows that $\kappa$ is measurable in an inner model [Fuchs].  Since not all universes have inner models with measurable cardinals, for some $V$ and $\kappa$ inaccessible, there is no extension where $\kappa$ is a singular worldly cardinal.  

The next theorem is about large cardinal properties just above inaccessible.  A cardinal $\kappa$ is $\Sigma_n$-\textit{correct} if and only if $V_{\kappa} \prec_{\Sigma_n} V$.  A cardinal $\kappa$ is $\Sigma_n$-\textit{reflecting} if and only if $\kappa$ is $\Sigma_n$-correct and inaccessible.  The following theorems show how force to change degrees of reflecting cardinals.  

\begin{theorem}  If $\kappa$ is $\Sigma_1$-reflecting, then there is a forcing extension where $\kappa$ is still $\Sigma_1$-reflecting, but not $\Sigma_2$-reflecting.  

\begin{proof}  First, let us establish that all inaccessible cardinals are $\Sigma_1$-reflecting.  Then, we shall prove the theorem by exhibiting a forcing notion which does not preserve that a cardinal is $\Sigma_2$-reflecting, but which does preserve that the cardinal is inaccessible (hence $\Sigma_1$-reflecting).  Suppose $\kappa$ is inaccessible. It follows that $V_{\kappa} = H_{\kappa}$.  The L\'{e}vy reflection theorem says that $H_{\kappa} \prec_{\Sigma_1} V$ for every uncountable cardinal $\kappa$.  Thus, if $\kappa$ is inaccessible, $V_{\kappa} \prec_{\Sigma_1} V$, and so it is $\Sigma_1$-reflecting.  Let $\mathbb{P}$ be the canonical forcing of the GCH, and let $G\subseteq \mathbb{P}$ be $V$-generic.  Next force over $V[G]$ to make $2^{\kappa} = \kappa^{++}$ in $V[G][g]$, where $g$ is $V[G]$-generic.  Since the canonical forcing of the GCH does not change that $\kappa$ is inaccessible, and forcing $2^{\kappa} = \kappa^{++}$ also cannot change that $\kappa$ is a strong limit (or regular), it follows that $\kappa$ is inaccessible in $V[G][g]$.  Thus $\kappa$ is $\Sigma_1$-reflecting in the final extension.  However, $V[G][g]_{\kappa} = V[G]_{\kappa} \models GCH$ while $V[G][g] \models 2^{\kappa} = \kappa^{++}$.  Thus, $\kappa$ is not $\Sigma_2$-reflecting in $V[G][g]$, since $V[G][g]_{\kappa}$ thinks the GCH holds everywhere while, $V[G][g]$ sees a violation of GCH at $\kappa$.  
\end{proof}
\end{theorem}

\newpage

\section{\textbf{Weakly compact cardinals and beyond}}

In this section, forcing notions are found which distinguish pairs of large cardinals in the [Mahlo, measurable] interval.  The first theorem separates the notions of Mahlo and weakly compact for a given cardinal with these properties.  A cardinal $\kappa$ is \textit{weakly compact} if and only if it is uncountable and satisfies the partition party $\kappa \to (\kappa)^2$. That is, $\kappa$ is weakly compact if and only if for every coloring of pairs of elements of $\kappa$ into two colors, there is a homogeneous subset of $\kappa$ of order-type $\kappa$.  There are other various equivalent characterizations of weakly compact cardinals. An important characterization is that weakly compact cardinals are inaccessible and have the tree property.  An uncountable cardinal $\kappa$ has the \textit{tree property} if every tree of height $\kappa$ whose levels have cardinality less than $\kappa$ has a branch of length $\kappa$ [Jech p. 120].  For example, $\omega_1$ does not have the tree property since there exists an Aronszajn tree, a tree of height $\omega_1$ with no uncountable branches.  A $\kappa$-Aronszajn tree is  a tree of height $\kappa$ with countable levels and no branches of length $\kappa$.  And, CH implies there is an $\omega_2$-Aranszajn tree, so in this universe, $\omega_2$ does not have the tree property.  A $\kappa$-Souslin tree is a tree of height $\kappa$ with no chains or antichains of size $\kappa$.  Thus a $\kappa$-Souslin tree is a $\kappa$-Aronszajn tree, since a chain is a branch.  Thus, the existence of a $\kappa$-Souslin tree destroys the tree property of $\kappa$, thus its weak compactness.

 The following theorems are analogues of Theorem \ref{Theorem. Mahlo every inaccessible degree.} and Theorem \ref{Theorem. Distinguish Mahlo and all inaccessible degrees.} from the previous section.  First, we see that every weakly compact cardinal is every degree of Mahlo.

\begin{theorem} If $\kappa$ is weakly compact then for every meta-ordinal term $t$ having only ordinals less than $\kappa$, the cardinal $\kappa$ is $t$-Mahlo.

\begin{proof}
Suppose $\kappa$ is weakly compact.  By induction suppose $\kappa$ is $s$-Mahlo for all $s < t$.  We want to show that $\kappa$ is $t$-Mahlo.  Fix any club $C \subseteq \kappa$.  Put $C \in M \models$ ZFC where $M$ is a $\kappa$-model, $M^{{<}\kappa} \subseteq M$.  Get $j:M \to N$ with $cp(j) = \kappa$.  Without loss of generality, we can assume $V_{\kappa} \subseteq M$.  Since degrees of Mahloness are downward absolute the fact that $\kappa$ is $s$-Mahlo in $V$ implies $\kappa$ is $s$-Mahlo in $N$.  Since $C$ is club, it follows that $\kappa \in j(C)$ so that $N \models$ ``the set of $s$-Mahlo cardinals meets $j(C)$''.  By elementarity, $M \models$ ``$C$ has an $s$-Mahlo cardinal''.  Thus $C$ has an $s$-Mahlo cardinal since $M$ is a $\kappa$-model.  Thus $\kappa$ is $t$-Mahlo.  
\end{proof}

\end{theorem}

Since the degrees of Mahlo cardinals discussed in the previous section are characterized by stationary sets, it follows that if a notion of forcing preserves stationary sets, it preserves any degree of Mahlo.  This theorem shows how to force a weakly compact cardinal to lose this property in a forcing extension while retaining that it is every possible degree of Mahlo.

\begin{theorem}
If $\kappa$ is weakly compact, then there is a forcing extension where for every meta-ordinal $t$ having only ordinals less than $\kappa$, the cardinal $\kappa$ is $t$-Mahlo but not weakly compact.
\begin{proof}
Suppose $\kappa$ is weakly compact. Let $\mathbb{S}$ be the forcing to add a $\kappa$-Souslin tree. Conditions in $\mathbb{S}$ are normal trees of height $\alpha +1$ for some $\alpha < \kappa$.  The conditions in $\mathbb{S}$ are ordered by end-extension: $s \le t$ if and only if $s$ is a tree which contains $t$ as an initial segment.  Also, it is required that for all $T \in \mathbb{S}$, the condition $T$ is a normal tree: every element of $T$ is a $\beta$-length binary sequence for some $\beta \le \alpha$, and if $t \subset s$ and $s \in T$ then $t \in T$, and every node splits in two at every level, and for every node $t$ and every level $\delta$ there is a node $s$ at level $\delta$ which extends $t$, and finally $|T| < \kappa.$  
Let $G \subset \mathbb{S}$ be $V$-generic. Let $S = \cup G$ be the new tree.  The generic object $S$ is a $\kappa$-Souslin tree [generalization of Jech p. 239].  These are the conditions.  The forcing notion $\mathbb{S}$ is $\kappa$-strategically closed since we can play the game along the paths as not to get stuck at an Aronszajn tree of height less than $\kappa$.  That is, the second player keeps track of the paths which are played and always chooses nodes which continue paths already chosen, so that at limits player two can continue play. Thus, $\mathbb{S}$ preserves stationary sets. Thus, the forcing preserves that $\kappa$ is any degree of Mahlo since the set of inaccessible cardinals, set of Mahlo cardinals, set of 1-Mahlo cardinals, etc, below $\kappa$, remain stationary in $\kappa$.  However, since the new tree, $S$, is a $\kappa$-Aronszajn tree, it does not have a branch of length $\kappa$.  Therefore, $\kappa$ no longer has the tree property, so it is no longer weakly compact.  

\end{proof}

\end{theorem}

A cardinal $\kappa$ is \textit{measurable} if and only if there is a transitive class $M$, and an elementary embedding $j:V \to M$ with critical point $\kappa$.  A cardinal $\kappa$ is \textit{weakly measurable} if and only if for every family $A \subseteq P(\kappa)$, of size at most $\kappa^+$, there is a non-principal $\kappa$-complete filter on $\kappa$ measuring every set in $A$. And, if $(\kappa)^{<\kappa} = \kappa$, then a cardinal $\kappa$ is weakly measurable if and only if for every transitive set $M$ of size $\kappa^+$, with $\kappa \in M$, there exists a transitive $N$ and an elementary embedding $j:M \to N$ with critical point $\kappa$. For the next large cardinal definition, we need the definition of a $\kappa$-model:  a transitive set $M$ is a $\kappa$-\textit{model} if and only if $M \models ZFC^{-}$ where $M$ is size $\kappa$ with $\kappa \in M$ and $M^{{<}\kappa} \subseteq M$.  Then, a cardinal $\kappa$ is \textit{strongly Ramsey} if every $A \subseteq \kappa$ is contained in a $\kappa$-model, $M$, for which there exists a weakly amenable $M$-ultrafilter on $\kappa$.  An $M$-ultrafilter $U$ is said to be weakly amenable (to $M$) if for every $A \in M$ of size $\kappa$ in $M$, the intersection $A \cap U$ is an element of $M$. An uncountable regular cardinal $\kappa$ is \textit{ineffable} if for every sequence $\langle A_{\alpha} \ | \ \alpha < \kappa \rangle$, with $A_{\alpha} \subseteq \alpha$, there is $A \subseteq \kappa$ such that the set $S = \{ \alpha < \kappa \ | \ A \cap \kappa = A_{\alpha} \}$ is stationary.  Measurable cardinals are weakly measurable; weakly measurable cardinals are strongly Ramsey cardinals; strongly Ramsey cardinals are not necessarily ineffable cardinals, but are below strongly Ramsey cardinals in consistency strength; and ineffable cardinals are weakly compact.  The following theorems show how to force the difference between these large cardinal properties for a given cardinal.  

The following theorem, by Jason Schanker [Schanker], fits exactly in the theme:

\begin{theorem} (Schanker) If $\kappa$ is weakly measurable, then the measurability of $\kappa$ can be destroyed while preserving that it is weakly measurable.
\end{theorem}

The preservation part of the following theorem, and the forcing used, is due to Victoria Gitman:

\begin{theorem} If $\kappa$ is strongly Ramsey, then there is a forcing extension where $\kappa$ is not weakly measurable, but is still strongly Ramsey.  

\begin{proof} Let $\kappa$ be strongly Ramsey. Make a preparation $\mathbb{P}$ to make $\kappa$ indestructible by adding a slim-$\kappa$-Kurepa tree [Gitman].  Let $G \subseteq \mathbb{P}$ be $V$-generic.  Let $\mathbb{K}$ be the forcing to add a slim-$\kappa$-Kurepa tree over $V[G]$. Let $H \subseteq \mathbb{K}$ be $V[G]$-generic.  In the extension, $\kappa$ is still strongly Ramsey [Gitman].  But, $\kappa$ cannot be weakly measurable in $V[G][H]$.  Suppose, by way of contradiction, that $\kappa$ is weakly measurable in the final extension.  Let $T \in V[G][H]$ be a slim-$\kappa$-Kurepa tree, and let $|M| = \kappa^+$, be a transitive set containing $\kappa$, the tree $T$, and the $\kappa^+$ many branches of $T$.  Since $\kappa$ is weakly measurable, there is an $N$, and an elementary embedding $j:M \to N$ with $cp(j) = \kappa$. Then, by elementarity, $j(T)$ is a slim-$j(\kappa)$-Kurepa tree.  But, since $T$ has $\kappa^+$ many paths, and $j(T)$ is an end-extension of $T$, each of the $\kappa^+$ many paths of $T$ extends to a node on the $\kappa$th level of $j(T)$.  Thus, there are $\kappa^+$ many nodes on the $\kappa$th level of $j(T)$, and $M$ can see all the paths of $T$.  So, $N$ can see that $j(T)$ is not a slim-$j(\kappa)$-Kurepa tree, which contradicts that $j(T)$ is a slim-$j(\kappa)$-Kurepa tree in $N$ by elementarity.  Thus, $\kappa$ is not weakly measurable in $V[G][H]$.
\end{proof}

\end{theorem}

The following theorem shows how to force to remove a large cardinal property while preserving is greater large cardinal property.  If a cardinal $\kappa$ is strongly Ramsey, then it is consistent that the cardinal is also an ineffable cardinal.  However, the following theorem shows how to force to an extension where $\kappa$ is not ineffable, but still strongly Ramsey.  

\begin{theorem}  If $\kappa$ is strongly Ramsey and ineffable, then there is a forcing extension where $\kappa$ is not ineffable, but it is still strongly Ramsey.

\begin{proof} Suppose $\kappa$ is ineffable.  The forcing to add a slim-$\kappa$-Kurepa tree, as in the previous theorem, will preserve that $\kappa$ is strongly Ramsey [Gitman].  However, forcing to add such a tree ensures that $\kappa$ is not ineffable in the extension since the characterization of $\kappa$ ineffable implies there are no slim-$\kappa$-Kurepa trees [Jensen/Kunen].
\end{proof}

\end{theorem}

\newpage
\section{\textbf{Measurable cardinals}}

In this section, degrees of measurable cardinals are changed by forcing. A cardinal $\kappa$ is \textit{measurable} if and only if $\kappa$ is uncountable and there is a $\kappa$-complete nonprincipal ultrafilter on $\kappa$.  Suppose $\kappa$ is measurable by measure $\mu$. Then, the corresponding ultrapower embedding $j: V \to M_{\mu}$, by measure $\mu$, is an elementary embedding with critical point $\kappa$.  Also, if $\kappa$ is the critical point of a non-trivial elementary embedding, then it is measurable.  Thus, another characterization of measurable cardinals is the following embedding definition: $\kappa$ is measurable if and only if it is the critical point of an elementary embedding $j: V \to M$ in $V$.  

The theorems in this section are about reducing Mitchell degrees of measurable cardinals in a forcing extension.  The Mitchell order is defined on normal measures; for normal measures $\mu$ and $\nu$, the following relation $u \lhd \nu$ holds when $\mu \in M_{\nu}$, where $M_{\nu}$ is the ultrapower by $\nu$.  If $\kappa$ is measurable, then the Mitchell order on the measures on $\kappa$ is well-founded (and transitive) [Jech p. 358]. Thus, let $m(\kappa)$ be the collection of normal measures on $\kappa$. The definition of rank for $\mu \in m(\kappa)$ with respect to $\lhd$ is $o(\mu) = \sup \{ o(\nu) + 1 \ | \nu \lhd \mu \}$.
  The Mitchell order on $\kappa$, denoted $o(\kappa)$, is the height of the Mitchell order on $m(\kappa)$, thus $o(\kappa) = \sup \{ o(\mu) + 1 \ | \ \mu \in m(\kappa) \}$.  If $\kappa$ is not measurable, then $o(\kappa) = 0$, since $m(\kappa)$ is empty.  A cardinal $\kappa$ is measurable if and only if $o(\kappa) \ge 1$. For $o(\kappa) \ge 2$, there is a measure $\mu$ on $\kappa$, such that if $j:V \to M_{\mu}$ is the corresponding ultrapower embedding, there is a measure on $\kappa$ in $M_{\mu}$.  Thus $\kappa$ is measurable in $M_{\mu}$.  Thus $M_{\mu} \models \kappa \in j(\{ \gamma < \kappa \ | \ \gamma \mbox{ is measurable } \} )$, and so the set $\{ \gamma < \kappa \ | \ \gamma \mbox{ is measurable } \}$ is in $\mu$.  That is, $\mu$ concentrates on measurable cardinals.  Therefore, $o(\kappa) \ge 2$ if and only if there is a normal measure on $\kappa$ which concentrates on measurable cardinals.  A cardinal $\kappa$ has Mitchell order $o(\kappa) \ge 3$, if and only if there is a normal measure on $\kappa$ which concentrates on measurable cardinals with Mitchell order 2.  A cardinal $\kappa$ has $o(\kappa) \ge \alpha$ if and only if for every $\beta < \alpha$, there is a normal measure on $\kappa$ which concentrates on measurable cardinals of order $\beta$.  

The following theorem shows how to force any measurable cardinal, of any Mitchell rank, to have order exactly one. In the proof of the following theorem, and the proof of the main theorems, we will need that the forcing notions do not create large cardinals.  In the previous sections, since being Mahlo or inaccessible is downwards absolute, we knew that the forcing could not create new Mahlo or inaccessible cardinals.  However, it is possible that forcing create large cardinals in the extension.  Therefore, we use the approximation and cover properties of [Hamkins4] to ensure that no new measurable cardinals are created.  

\begin{theorem} If $\kappa$ is a measurable cardinal, then there is a forcing extension where $o(\kappa) = 1$.  In addition, the forcing will preserve all measurable cardinals and will not create any new measurable cardinals.

\begin{proof} 
Suppose $\kappa$ is a measurable cardinal.  Assume that $2^{\kappa} = \kappa^+$, or force it to be true since forcing to have the GCH hold at $\kappa$ preserves the measurability of $\kappa$.  Pick $j: V \to M$, an elementary ultrapower embedding by $\mu$, a normal measure on $\kappa$ which has minimal Mitchell rank $o(\mu) = 0$, with $cp(j) = \kappa$, such that $M \models \kappa$ is not measurable.  Let $\mathbb{P}$ be a $\kappa$-length iteration, with Easton support, of $\mathbb{Q}_{\gamma}$, a $\mathbb{P}_{\gamma}$-name for the forcing to add a club to $\gamma$ of cardinals which are not measurable in $V$, whenever $\gamma < \kappa$ is inaccessible, otherwise force trivially at stage $\gamma$.  Let $G \subseteq \mathbb{P}$ be $V$-generic, let $\dot{\mathbb{Q}}_G = \mathbb{Q}$, and let $\dot{\mathbb{Q}}$ be a $\mathbb{P}$-name for a forcing to add a club of cardinals to $\kappa$ which are not measurable in $V[G]$. Let $g \subseteq \mathbb{Q}$ be $V[G]$-generic.  Let $\delta_0$ is the first inaccessible cardinal.  This is the first non-trivial stage of $\mathbb{P}$, and at this stage a club is added to $\delta_0$.  Thus the forcing up to and including stage $\delta_0$ has cardinality at most $\delta_0$ and the forcing after stage $\delta_0$ is closed up to the next inaccessible cardinal.  Thus, $\mathbb{P}$ has a closure point at $\delta_0$.  Thus by Corollary 22 in [Hamkins4], the forcing does not create measurable cardinals.

We want to lift $j$ through $\mathbb{P}$ and $\mathbb{Q}$ to find $j: V[G][g] \to M[j(G)][j(g)]$ which witnesses that $\kappa$ is measurable in $V[G][g]$. The main task is to find $j(G)$, an $M$-generic filter for $j(\mathbb{P})$, and $j(g)$ an $M[j(G)]$-generic filter for $j(\mathbb{Q})$.  Since the critical point of $j$ is $\kappa$, the forcings $\mathbb{P}$ and $j(\mathbb{P})$ are isomorphic up to stage $\kappa$. Also, since $\dot{\mathbb{Q}}$ is a $\mathbb{P}$-name for a forcing which adds a club of ground model non-measurable cardinals to $\kappa$ and $V[G]$ and $M[G]$ have the same bounded subsets of $\kappa$, the forcing at stage $\kappa$ of $j(\mathbb{P})$ is $\dot{\mathbb{Q}}$.  Thus, $j(\mathbb{P})$ factors as $\mathbb{P} \ast \dot{\mathbb{Q}} \ast \mathbb{P}_{\text{tail}}$, where $\mathbb{P}_{\text{tail}}$ is the forcing $j(\mathbb{P})$ past stage $\kappa$.  Since $G$ is $V$-generic, it is also $M$-generic.  Thus, form the structure $M[G]$.  Similarly, since $g$ is $V[G]$-generic, it follows that $g$ is $M[G]$-generic, so construct the structure $M[G][g]$.  Since both $\mathbb{P}$ and $\mathbb{Q}$ have size $\kappa$, they both have the $\kappa^+$-chain condition, so the iteration $\mathbb{P} \ast \mathbb{Q}$ has the $\kappa^+$-chain condition.  Also, since $j''\kappa = \kappa \in M$, and $j$ is an ultrapower embedding, it follows that $M^{\kappa} \subseteq M$ in $V$.  Thus, since $G \ast g$ is $V$-generic, $M[G][g]^{\kappa} \subseteq M[G][g]$ [Hamkins1 (Lemma 54)].

The goal is to diagonalize to find an $M[G][g]$-generic filter for $\mathbb{P}_{\text{tail}}$ in $V[G][g]$, and so far we have satisfied one criterion which is that $M[G][g]$ is closed under $\kappa$ sequences in $V[G][g]$.  Next, since $|\mathbb{P}| = \kappa$, and $2^{\kappa} = \kappa^+$, it follows that $\mathbb{P}$ has at most $\kappa^+$ dense sets, and so does any tail $\mathbb{P}_{[\beta, \kappa)}$.  Thus, $\mathbb{P}_{\text{tail}}$ has at most $|j(\kappa^+)|^V \le \kappa^{+^{\kappa}} = \kappa^+$ dense subsets in $M[G][g]$.  Finally, since for every $\beta < \kappa$, there is, in $M[G][g]$, a dense subset of $\mathbb{P}$ which is $\le \beta$-closed, there is a dense subset of $\mathbb{P}_{\text{tail}}$ which is $\le \kappa$-closed. Thus, it is possible to diagonalize to get an $M[G][g]$-generic filter $G_{\text{tail}} \subseteq \mathbb{P}_{\text{tail}}$, in $V[G][g]$.  Thus, $j(G) = G \ast g \ast G_{\text{tail}} \subseteq j(\mathbb{P})$ is $M$-generic, and $j''G \subseteq G \ast g \ast G_{\text{tail}}$.  Therefore, the lifting criterion is satisfied, and so the embedding lifts to $j:V[G] \to M[j(G)]$.  

The next goal is to lift the embedding through $\mathbb{Q}$.  The forcing $j(\mathbb{Q})$ is the forcing to add a club of non-measurable cardinals to $j(\kappa)$, and we need an $M[j(G)]$-generic filter for $j(\mathbb{Q})$.  Since $G_{\text{tail}} \in V[G][g]$, it follows that $M[j(G)]^{\kappa} \subseteq M[j(G)]$.  Also, since $|\mathbb{Q}| = \kappa$, it has at most $\kappa^+$ many dense sets in $V[G]$, by elementarity $j(\mathbb{Q})$ has at most $|j(\kappa^+)|^{V[G]} \le \kappa^{+^{\kappa}} = \kappa^+$ many dense subsets in $M[j(G)]$.  For any $\beta < \kappa$, the forcing $\mathbb{Q}$ has a dense subset which is $\le \beta$-closed. Since $\kappa < j(\kappa)$, the forcing $j(\mathbb{Q})$ has a dense subset which is $\le \kappa$-closed.  Note that $c = \cup g$, which is club in $\kappa$, is in $M[j(G)]$. Let $\bar{c} = c \cup \{ \kappa \}$.  For all $\delta < \kappa$, we have $\delta \in c$ then $\delta$ is not measurable in $V$.  Thus, $\delta \in c$ implies $\delta$ is not measurable in $M$, since $M$ and $V$ have the same $P(P(\delta))$ since $\delta < \kappa$.  By Corollary 22 in [Hamkins4] since $M \subseteq M[j(G)]$ also satisfies approximation and cover properties, it follows that $\delta$ is not measurable in $M[j(G)]$.  By choice of $j:V \to M$, we have that $\kappa$ is not measurable in $M$.  Thus $\kappa$ is not measurable in $M[j(G)]$, again by [Hamkins4]. Thus $\bar{c}$ is a closed, bounded subset of $j(\kappa)$ which contains no measurable cardinals.  Thus $\bar{c} \in j(\mathbb{Q})$. Thus, diagonalize to get a generic filter $g^{\ast} \subseteq j(\mathbb{Q})$, in $V[G]$, which meets a dense subset which is $\le \kappa$-closed, and which contains $c$.  Thus, let $j(g) = g^{\ast}$, and observe that $j''g \subseteq g^{\ast}$ since $g \in \cup g^{\ast}$.  Thus, the embedding lifts to $j: V[G][g] \to M[j(G)][j(g)]$, in $V[G][g]$.  Therefore, $\kappa$ is still measurable in $V[G][g]$.  Thus $o(\kappa) \ge 1$ in $V[G][g]$.

Finally $o(\kappa)^{V[G][g]} \le 1$ will follow from the fact that $\mathbb{Q}$ adds a club subset $c \subseteq \kappa$ which contains no cardinals which are measurable in $V$. Since the combined forcing does not create measurable cardinals, the club $c$ contains no cardinals which are measurable in $V[G][g]$.  Since the new club has measure one in any normal measure on $\kappa$, the complement of $c$, which contains the measurable cardinals of $V[G][g]$ below $\kappa$ has measure zero.  Thus, there is no normal measure on $\kappa$ in $V[G][g]$ which concentrates on the measurable cardinals of $V[G][g]$.  Hence $o(\kappa)^{V[G][g]} \ngtr 1.$  Thus, $o(\kappa) = 1$ in $V[G][g]$.  That the forcing preserves all measurable cardinals and creates no new measurable cardinals will follow from the proof of the more general Theorem \ref{Theorem. Kill measurable softly.}.
\end{proof}
\end{theorem}

The previous theorem shows how to make the Mitchell order of any measurable cardinal exactly one, for any measurable cardinal.  The next theorem will show how to make the Mitchell order of any measurable cardinal exactly $\alpha$, for any $\alpha < \kappa$, in a forcing extension, whenever the order is at least $\alpha$. Preceding the proof are the following lemmas, so that we may generalize the proof for any $\alpha < \kappa$. 

\begin{lemma} \label{Lemma. Club shooting.}  Let $\kappa$ be an infinite cardinal and $\kappa^{<\kappa} = \kappa$.  Let $S \subseteq \kappa$ be a subset of $\kappa$ which contains the singular cardinals.  Then $\mathbb{Q}_S$, the forcing to add a club $C \subset S$, preserves cardinals and cofinalities, and for all $\beta < \kappa$, the forcing $\mathbb{Q}_S$ has a $\le \beta$-closed dense subset.

\begin{proof} Conditions $c \in \mathbb{Q}_S$ are closed, bounded subsets of $S$.  The ordering on $\mathbb{Q}_S$ is end-extension: $c \le d$ if and only if $a \in c \setminus d \rightarrow a > \max(d)$.  Let $G \subseteq \mathbb{Q}_S$ be $V$-generic, and let $C = \cup G$.  Then $C \subseteq S$  is club in $\kappa$.  First, see that $C$ is unbounded.  Let $\beta \in \kappa$, and let $D_{\beta} = \{ d \in \mathbb{Q}_S \ : \ \max(d) > \beta \}$.  Then, $D_{\beta}$ is dense in $\mathbb{Q}_S$, since if $c \in \mathbb{Q}_S$, the following is a condition: $ d = c \cup \{ \beta' \} $, where $\beta' > \beta$ and $\beta'$ is in $S$.  Then, $d$ is an end-extension of $c$, hence stronger, and $d \in D_{\beta}$, which shows $D_{\beta}$ is dense in $\mathbb{Q}_S$.  Hence, there is $c_{\beta} \in D_{\beta} \cap G$, so that there is an element of $C$ above $\beta$.  Thus, $C$ is unbounded in $\kappa$.  Next, see that $C$ is closed.  Suppose $\delta \cap C$ is unbounded in $\delta < \kappa$.  Since $C$ is unbounded, there is $\delta' \in C$, where $\delta' > \delta$, and a condition $c \in C$ which contains $\delta'$.  Since the conditions are ordered by end-extension, and since there is a sequence (possibly of length 1) of conditions which witness that $\delta \cap C$ is unbounded in $C$, which will all be contained in the condition $c$, which contains an element above $\delta$, it follows that $\delta \cap c$ is unbounded in $\delta$.  Since $c$ is closed, $\delta \in c$.  Therefore $\delta \in C$, which shows $C$ is closed.  

Furthermore, $\mathbb{Q}_S$ preserves cardinals and cofinalities $\ge \kappa^+$ since $|\mathbb{Q}_S| = \kappa^{<\kappa} = \kappa$.  If $\beta < \kappa$, then $\mathbb{Q}_S$ is not necessarily $\le \beta$-closed.  Let $\beta$ be the first limit of $S$ which is not in $S$.  Then, the initial segments of $S$ below $\beta$ are in $\mathbb{Q}_S$, and union up to $S \cap \beta$, so that no condition can close this $\beta$-sequence since $\beta \notin S$.  However, for every $\beta < \kappa$, the set $D_{\beta} = \{ d \in \mathbb{Q}_S \ : \ \max(d) \ge \beta \}$ is dense in $\mathbb{Q}_S$ and is $\le \beta$-closed since cof($\kappa) > \beta$.  Thus, for every $\beta < \kappa$, the poset $\mathbb{Q}_S$ is forcing equivalent to $\mathbb{C}_{\beta}$, which is $\le \beta$-closed.  Thus, $\mathbb{Q}_S$ preserves all cardinals, cofinalities, and strong limits $\le \kappa$.

\end{proof}
\end{lemma}

The following Lemma is about the Mitchell rank of a normal measure $\mu$ on a measurable cardinal $\kappa$. By definition, $o(\kappa)$ is the height of the well-founded Mitchell relation $\lhd$ on $m(\kappa)$, the collection of normal measures on $\kappa$. Thus $o(\kappa) = \{ o(\mu) + 1 \ | \ \mu \in m(\kappa) \}$.  Similarly the definition of rank for $\mu \in m(\kappa)$ with respect to $\lhd$ is $o(\mu) = \sup \{ o(\nu) + 1 \ | \nu \lhd \mu \}$.

\begin{lemma} \label{Lemma. Mitchell rank in M.}  If $\kappa$ is a measurable cardinal, and $\mu$ is a normal measure on $\kappa$, and $j_{\mu} : V \to M_{\mu}$ is the ultrapower embedding by $\mu$ with critical point $\kappa$, then $o(\mu) = o(\kappa)^{M_{\mu}}$.  Thus, $\beta < o(\kappa)$ if and only if there exists $j:V \to M$ elementary embedding with critical point $\kappa$ with $o(\kappa)^M = \beta$.

\begin{proof}
Suppose $\kappa$ is measurable, $\mu$ is a normal ulltrafilter on $\kappa$, and $j_{\mu}: V \to M_{\mu}$ is the ultrapower embedding by $\mu$ with critical point $\kappa$.   The model $M_{\mu}$ computes $\lhd$ correctly for $\nu \lhd \mu$, which means $\nu \in M_{\mu}$, since $M^{\kappa} \subseteq M$ and so $M_{\mu}$ has all the functions $f:\kappa \to V_{\kappa}$ needed to check whether $\nu \lhd \mu$.  Then
\begin{align*}
    o(\kappa)^{M_{\mu}}  &= \sup\{ o(\nu)^{M_{\nu}} + 1 \ | \ \nu \in M_{\mu} \} \\ 
     &= \sup\{o(\nu)^{M_{\mu}} + 1 \ | \ \nu \lhd \mu \} \\
     &= \sup\{o(\nu) + 1 \ | \ \nu \lhd \mu \} \\
     &= o(\mu)
\end{align*}
Thus, the Mitchell rank of $\kappa$ in $M_{\mu}$ is the Mitchell rank of $\mu$.  Let $\beta < o(\kappa)$.  Since $\lhd$ is well-founded, there is $\mu$ a normal measure on $\kappa$ such that $o(\mu) = \beta$.  Then if $j_{\mu} : V \to M_{\mu}$ is the ultrapower by $\mu$ with critical point $\kappa$, then $o(\kappa)^{M_{\mu}} = o(\mu) = \beta$, as desired.  Conversely, suppose there is an elementary embedding $j:V \to M$ with cp$(j) = \kappa$ and $o(\kappa)^M = \beta$.  Let $\mu$ be the induced normal measure, using $\kappa$ as a seed ($X \in \mu$ if and only if $\kappa \in j(X)$) and let $M_{\mu}$ be the ultrapower by $\mu$.  Then $M$ and $M_{\mu}$ have the same normal measures and the same Mitchell order and ranks, so $M_{\mu} \models o(\kappa) = \beta$.  Since we have established $o(\kappa)^{M_{\mu}} = o(\mu)$, the Mitchell rank of $\mu$ in $m(\kappa)$ is $\beta$, that is $o(\mu) = \beta$.  Thus since $o(\kappa) = \{ o(\mu) + 1 | \mu \in m(\kappa) \}$, it follows that $o(\kappa)^V > \beta$.

\end{proof}
\end{lemma}

The following lemma generalizes Corollary 22 of [Hamkins4] by showing that if $V \subseteq \bar{V}$ satisfies the $\delta$ approximation and cover properties, then the Mitchell rank of a cardinal $\kappa > \delta$ cannot increase between $V$ and $\bar{V}$.

\begin{lemma} \label{Lemma. Mitchell rank doesn't go up.} Suppose $V \subseteq \bar{V}$ satisfies the $\delta$ approximation and cover properties.  If $\kappa > \delta$, then $o(\kappa)^{\bar{V}} \le o(\kappa)^V$.  

\begin{proof}  Suppose $V \subseteq \bar{V}$ satisfies the $\delta$ approximation and cover properties.  Suppose $\kappa > \delta$.  We want to show $o(\kappa)^{\bar{V}} \le o(\kappa)^V$ by induction on $\kappa$.  Assume inductively that for all $\kappa' < \kappa$, we have $o(\kappa')^{\bar{V}} \le o(\kappa')^V$.  If $o(\kappa)^{\bar{V}} > o(\kappa)^V$ then by Lemma \ref{Lemma. Mitchell rank in M.} get and elementary ultrapower embedding $j: \bar{V} \to \bar{M}$ with critical point $\kappa$ and $o(\kappa)^{\bar{M}} = o(\kappa)^{V}$.  By [Hamkins4], $j$ is a lift of an ultrapower embedding $j \restriction V: V \to M$ for some $M$, and this embedding is a class in $V$. Again by \ref{Lemma. Mitchell rank in M.} we have established that $o(\kappa)^M < o(\kappa)^V$.  Thus, $o(\kappa)^M < o(\kappa)^{\bar{M}}$.  However, this contradicts that $j$ applied to the induction hypothesis gives $o(\kappa)^{\bar{M}} \le o(\kappa)^{M}$ since $\kappa < j(\kappa)$.  Therefore, $o(\kappa)^{\bar{V}} \le o(\kappa)^V$.
\end{proof}

\end{lemma}

\begin{lemma} \label{Lemma. Lift normal ultrapower embeddings.} If $V \subseteq V[G]$ is a forcing extension and every normal ultrapower embedding $j:V \to M$ in $V$ (with any critical point $\kappa$) lifts to an embedding $j:V[G] \to M[j(G)]$, then $o(\kappa)^{V[G]} \ge o(\kappa)^V$.  In addition, if $V \subseteq V[G]$ also has $\delta$ approximation and $\delta$ cover properties, then $o(\kappa)^{V[G]} = o(\kappa)^V$.
\begin{proof} Suppose $V \subseteq V[G]$ is a forcing extension where every normal ultrapower embedding in $V$ lifts to a normal ultrapower embedding in $V[G]$.  Suppose $\kappa$ is a measurable cardinal, and for every measurable cardinal $\kappa' < \kappa$, assume $o(\kappa')^{V[G]} \ge o(\kappa')^V$.  Fix any $\beta < o(\kappa)$.  By Lemma \ref{Lemma. Mitchell rank in M.} there exists an elementary embedding $j:V \to M$ with $cp(j) = \kappa$ and $M \models o(\kappa) = \beta$.  Then, by the assumption on $V \subseteq V[G]$, this $j$ lifts to $j:V[G] \to M[j(G)]$.  Then, applying $j$ to the induction hypothesis gives $o(\kappa)^{M[j(G)]} \ge o(\kappa)^M$ since $\kappa < j(\kappa)$.  Then, since $M \models o(\kappa) = \beta$ it follows that $o(\kappa)^{M[j(G)]} \ge \beta$.  Thus by Lemma \ref{Lemma. Mitchell rank in M.} $o(\kappa)^{V[G]} > \beta$.  Thus $o(\kappa)^{V[G]} \ge o(\kappa)^V$.  

If $V \subseteq V[G]$ also satisfies the $\delta$ approximation and $\delta$ cover properties, then by Lemma \ref{Lemma. Mitchell rank doesn't go up.} the Mitchell rank does not go up between these models. Thus $o(\kappa)^{V[G]} = o(\kappa)^V$.
\end{proof}

\end{lemma} 

\begin{theorem} \label{Theorem. Kill measurable softly.} For any $V \models ZFC + GCH$ and any ordinal $\alpha$, there is a forcing extension $V[G]$ where every cardinal $\kappa$ above $\alpha$ in $V[G]$ has $o(\kappa)^{V[G]} = \min \{ \alpha, o(\kappa)^V \}$.  In particular, if $1 \le \alpha$, then the forcing preserves all measurable cardinals of $V$ and creates no new measurable cardinals.

\begin{proof}
Let $\alpha \in Ord$.  The following notion of forcing will ensure that the Mitchell rank of all cardinals above $\alpha$ which have Mitchell rank above $\alpha$ will change. Let $\mathbb{P}$ be an Easton support Ord-length iteration, forcing at inaccessible stages $\gamma$, to add a club $c_{\gamma} \subseteq \gamma$ such that $\delta \in c_{\gamma}$ implies $o(\delta)^V < \alpha$. Let $G \subseteq \mathbb{P}$ be $V$-generic. Let $\delta_0$ denote the first inaccessible cardinal.  The forcing up to this stage, $P_{\delta_0}$, adds nothing and the forcing $\mathbb{Q}$ at stage $\delta_0$ adds a club to $\delta_0$, so that $|\mathbb{P}_{\delta_0} \ast \dot{\mathbb{Q}}| \le \delta_0$.  The forcing past stage $\delta_0$ is closed up to the next inaccessible cardinal.  Since $\mathbb{P}$ has a closure point at $\delta_0$, by Lemma 13 in [Hamkins4] it follows that $V \subseteq V[G]$ satisfies the $\delta_0^+$ approximation and cover properties.  Thus, by Lemma \ref{Lemma. Mitchell rank doesn't go up.} the Mitchell rank does not go up between $V$ and $V[G]$.  Note that any measurable cardinal is above $\delta_0^+$.  

Let $\kappa$ be a measurable cardinal above $\alpha$.  Fix any $\beta_0 < \min\{ \alpha, o(\kappa)^V \}$.  The induction hypothesis is that for all $\kappa' < \kappa$, we have $o(\kappa')^{V[G]} = \min\{ \alpha, o(\kappa')^V\}$.  By Lemma \ref{Lemma. Mitchell rank in M.}, fix an elementary embedding $j:V \to M$ with $cp(j) = \kappa$ and $o(\kappa)^M = \beta_0$.  We shall lift $j$ through $\mathbb{P}$.  The forcing above $\kappa$ does not affect $o(\kappa)$ since the forcing past stage $\kappa$ has a dense subset which is closed up to the next inaccessible.  Allow me to call $\mathbb{P}_{\kappa} \ast \dot{\mathbb{Q}}$ the forcing up to and including stage $\kappa$ while I am lifting the embedding, first through $\mathbb{P}_{\kappa}$ and then through $\mathbb{Q}$.  Let $G_{\kappa} \subseteq \mathbb{P}_{\kappa}$ be $V$-generic, and $g \subseteq \mathbb{Q}$ be $V[G_{\kappa}]$-generic.  First, we need an $M$-generic filter for $j(\mathbb{P}_{\kappa})$.  Since $cp(j) = \kappa$, the forcings $\mathbb{P}_{\kappa}$ and $j(\mathbb{P}_{\kappa})$ agree up to stage $\kappa$.  Since $\alpha < \kappa$, where $\alpha$ is the ordinal in the statement of the theorem, and the forcing $\mathbb{Q}$ adds a club $c_{\kappa} \subseteq \kappa$ such that $\delta \in c_{\kappa}$ implies $o(\delta)^V < \alpha$ which implies $o(\delta)^M < \alpha$, it follows that the $\kappa$th stage of $j(\mathbb{P}_{\kappa})$ is $\mathbb{Q}$.  Thus, $j(\mathbb{P}_{\kappa})$ factors as $\mathbb{P}_{\kappa} \ast \dot{\mathbb{Q}} \ast \mathbb{P}_{\text{tail}}$ where $\mathbb{P}_{\text{tail}}$ is the forcing past stage $\kappa$.  Since $G_{\kappa} \ast g$ is generic for $\mathbb{P}_{\kappa} \ast \dot{\mathbb{Q}}$, we just need an $M[G_{\kappa}][g]$-generic filter for $\mathbb{P}_{\text{tail}}$ which is in $V[G_{\kappa}][g]$ (since we want the final embedding to be there).  Thus, diagonalize to get an $M[G_{\kappa}][g]$-generic filter for $\mathbb{P}_{\text{tail}}$ after checking we have met the criteria.  Namely, since $|\mathbb{P}_{\kappa}| = |\dot{Q}| = \kappa$, the forcing $\mathbb{P}_{\kappa} \ast \dot{\mathbb{Q}}$ has the $\kappa^+$-chain condition.  Since $M^{\kappa} \subseteq M$, it follows [Hamkins1 (Theorem 54)] that $M[G_{\kappa}][g]^{\kappa} \subseteq M[G_{\kappa}][g]$.  Also, since $2^{\kappa} = \kappa^+$, the forcing $\mathbb{P}_{\kappa}$ has $\kappa^+$ many dense sets in $V$.  Thus, $\mathbb{P}_{\text{tail}}$ has at most $|j(\kappa^+)|^V \le \kappa^{+^{\kappa}} = \kappa^+$ many dense sets in $M[G_{\kappa}][g]$.  Also, since for any $\beta < \kappa$, the forcing $\mathbb{P}_{\kappa}$ has a dense subset which is ${\le}\beta$-closed (by Lemma \ref{Lemma. Club shooting.}), it follows that there is a dense subset of $\mathbb{P}_{\text{tail}}$ which is ${\le}\kappa$-closed.  Thus, we may diagonalize to get $G^* \subset \mathbb{P}_{\text{tail}}$ for this dense subset, a generic filter in $V[G_{\kappa}][g]$, which is $M[G_{\kappa}][g]$-generic.  Thus, $j(G_{\kappa}) = G_{\kappa} \ast g \ast G^*$ is $M$-generic for $j(\mathbb{P}_{\kappa})$ and $j''G_{\kappa} \subseteq G_{\kappa} \ast g \ast G^*$, so we may lift $j$ to $j:V[G_{\kappa}] \to M[j(G_{\kappa})]$.  

Next, we lift $j$ through $\mathbb{Q}$.  The forcing $j(\mathbb{Q})$ adds a club $C \subseteq j(\kappa)$ such that $\delta \in C$ implies $o(\delta)^{M[j(G_{\kappa})]} < \alpha$.  Since $G_{\kappa} \ast g \ast G^* \in V[G_{\kappa}][g]$ and $M^{\kappa} \subseteq M$, it follows [Hamkins1 (Theorem 53)] that $M[j(G_{\kappa})]^{\kappa} \subseteq M[j(G_{\kappa})]$.  Also, by Lemma \ref{Lemma. Club shooting.}, there is a dense subset of $j(\mathbb{Q})$ which is ${\le}\kappa$-closed, and since $|\mathbb{Q}| = \kappa$, it follows that $j(\mathbb{Q})$ has at most $\kappa^+$ many dense subsets in $M[j(G_{\kappa})]$.  Let $c = \cup g$, which is club in $\kappa$, and consider the set $\bar{c} = c \cup \{ \kappa \} \in M[j(G_{\kappa})]$.  We have $\delta \in c$ implies $o(\delta)^V < \alpha$. Since $M$ and $V$ agree on $P(P(\delta))$ for $\delta < \kappa$, it follows that $o(\delta)^V = o(\delta)^M$. By Lemma \ref{Lemma. Mitchell rank doesn't go up.}, we have $o(\delta)^{M[j(G_{\kappa})]} \le o(\delta)^M$.  Thus $\delta \in c$ implies $o(\delta)^{M[j(G_{\kappa})]} < \alpha$.  Since $o(\kappa)^M < \alpha$, it follows that $o(\kappa)^{M[j(G_{\kappa})]} < \alpha$, so that $\bar{c}$ is a condition of $j(\mathbb{Q})$, since $\bar{c}$ is a closed, bounded subset of $j(\kappa)$ which such that $\delta \in \bar{c}$ implies $o(\delta)^{M[j(G_{\kappa})]} < \alpha$. Thus, diagonalize to get an $M[j(G)]$-generic filter, $g^* \subseteq j(\mathbb{Q})$ which contains the condition $\bar{c}$. Then $j''g \subseteq j(g) = g^*$. Therefore, we may lift the embedding to $j: V[G_{\kappa}][g] \to M[j(G_{\kappa})][j(g)]$.

From the previous lifting arguments will follow that $o(\kappa)^{V[G]} \ge \min \{ o(\kappa)^V, \alpha \}$.  The lifted embedding $j:V \to M$ was chosen so that $o(\kappa)^M = \beta_0 < \min \{ o(\kappa)^V, \alpha \}$.  Since $\beta_0$ was arbitrary, it follows from Lemma \ref{Lemma. Mitchell rank in M.} that $o(\kappa)^V \ge \min \{ o(\kappa)^V, \alpha \}$.  All that remains to see is $o(\kappa)^{V[G]} \le \min \{ o(\kappa)^V, \alpha \}$.  By Lemma \ref{Lemma. Mitchell rank doesn't go up.}, which shows Mitchell rank does not go up for $V \subseteq V[G]$, we have $o(\kappa)^{V[G]} \le o(\kappa)^V$.  Consider any embedding $j:V[G] \to M[j(G)]$ with critical point $\kappa$ in $V[G]$, and the new club $c \subseteq \kappa$, where $\forall \delta < \kappa$, $\delta \in c$ implies $o(\delta)^V < \alpha$, which is in any normal measure on $\kappa$.  By the induction hypothesis, $\delta \in c$ implies $o(\delta)^{V[G]} < \alpha$.  Applying $j$ to this statement gives, $\forall \delta < j(\kappa)$, $\delta \in j(c)$ implies $o(\delta)^{M[j(G)]} < \alpha$.  Since $\kappa < j(\kappa)$ and $\kappa \in j(c)$, this gives $o(\kappa)^{M[j(G)]} < \alpha$.  Thus by Lemma \ref{Lemma. Mitchell rank in M.}, since $j$ was arbitrary, $o(\kappa)^{V[G]} \le \alpha$.   Thus $o(\kappa)^{V[G]} \le \min \{ o(\kappa)^V, \alpha \}$.  Thus, $o(\kappa)^{V[G]} = \min \{ o(\kappa)^V, \alpha \}$.
\end{proof}
\end{theorem}

I would like to now consider cases where $o(\kappa) \ge \kappa$. For this we use the concept of representing function.  Suppose $\kappa$ is measurable, $o(\kappa) = \alpha$, and $\alpha \in H_{\theta^+}$.  
The ordinal $\alpha$ is represented by $f: \kappa \to V_{\kappa}$ if whenever $j: V \to M$ an elementary embedding with cp$(j) = \kappa$ and $M^{\kappa} \subseteq M$, then $j(f)(\kappa) = \alpha$ (this can be formalized as a first-order statement using extenders) [Hamkins2].  Let $F: \text{Ord} \to \text{Ord}$. The main theorem shows how to change the Mitchell rank for all measurable cardinals $\gamma$ for which $F \restriction \gamma$ represents $F(\gamma)$ (through any $\gamma$-closed embedding with critical point $\gamma$). 

\begin{theorem} For any $V \models ZFC + GCH$ and any $F: \text{Ord} \to \text{Ord}$, there is a forcing extension $V[G]$ where, if $\kappa$ is measurable and $F \restriction \kappa$ represents $F(\kappa)$ in V, then $o(\kappa)^{V[G]} =\min\{o(\kappa)^V, F(\kappa)\}$.
\begin{proof}
For example, we can do this for numerous functions satisfying this hypothesis, such as $F(\gamma) = \gamma$, $F(\gamma) = \gamma + 3$, $F(\gamma) = \gamma^2 + \omega^2 \cdot 2 + 5$, and any other function which can be defined like this.  Suppose $V \models \text{ZFC} + \text{GCH}$. Let $F : \text{Ord} \to \text{Ord}$. Let $\mathbb{P}$ be an Easton support Ord-length iteration, forcing at inaccessible stages $\gamma$, to add a club $C_{\gamma} \subseteq \gamma$ such that $\delta \in C_{\gamma}$ implies $o(\delta)^V < F(\delta)$.  Let $\delta_0$ denote the first inaccessible cardinal.  The forcing up to this stage, $P_{\delta_0}$, adds nothing and the forcing $\mathbb{Q}$ at stage $\delta_0$ adds a club to $\delta_0$, so that $|\mathbb{P}_{\delta_0} \ast \dot{\mathbb{Q}}| \le \delta_0$.  The forcing past stage $\delta_0$ is closed up to the next inaccessible cardinal.  Let $G \subseteq \mathbb{P}$ be $V$-generic.  Since $\mathbb{P}$ has a closure point at $\delta_0$, by Lemma 13 in [Hamkins4] it follows that $V \subseteq V[G]$ satisfies the $\delta_0^+$ approximation and cover properties.  Thus, by Lemma \ref{Lemma. Mitchell rank doesn't go up.} the Mitchell rank does not go up between $V$ and $V[G]$.  

Let $\kappa$ be a measurable cardinal such that $F \restriction \kappa$ represents $F(\kappa) = \alpha$ in $V$. Here we will need that representing functions are still representing functions in the extension.  The fact we will use for this is from [Hamkins4] which says that if $V \subseteq V[G]$ satisfies the hypothesis of approximation and cover theorem, then every $\kappa$-closed embedding $j:V[G] \to M[j(G)]$ is a lift of an embedding $j \restriction V : V \to M$ that is in $V$ (and $j \restriction V$ is also a $\kappa$-closed embedding.  Thus, if $F \restriction \kappa$ is a representing function for $F(\kappa)$ in the ground model, $(j \restriction V) (F \restriction \kappa)(\kappa) = F(\kappa) = \alpha$, and thus $F \restriction \kappa$ is a representing function for $F(\kappa)$ with respect to such embeddings $j$ in $V[G]$.

The induction hypothesis is that for all $\kappa' < \kappa$ for which $F \restriction \kappa'$ represents $F(\kappa')$ in V has $o(\kappa')^{V[G]} = \min\{ F(\kappa'), o(\kappa')^V\}$.  By Lemma \ref{Lemma. Mitchell rank in M.} fix an elementary embedding $j:V \to M$ with $cp(j) = \kappa$ and $o(\kappa)^M = \beta_0 < \min \{o(\kappa)^V, F(\kappa) \}$. Note that $j(F)(\kappa) = j(F \restriction \kappa)(\kappa) = F(\kappa)$.  We shall lift $j$ through $\mathbb{P}$.  The forcing above $\kappa$ does not affect $o(\kappa)$ since the forcing past stage $\kappa$ is closed up to the next inaccessible.  Allow me to call $\mathbb{P}_{\kappa} \ast \dot{\mathbb{Q}}$ the forcing up to and including stage $\kappa$ while I am lifting the embedding, first through $\mathbb{P}_{\kappa}$ and then through $\mathbb{Q}$.  Let $G_{\kappa} \subseteq \mathbb{P}_{\kappa}$ be $V$-generic, and $g \subseteq \mathbb{Q}$ be $V[G_{\kappa}]$-generic.  First, we need an $M$-generic filter for $j(\mathbb{P}_{\kappa})$.  Since $cp(j) = \kappa$, the forcings $\mathbb{P}_{\kappa}$ and $j(\mathbb{P}_{\kappa})$ agree up to stage $\kappa$.  The forcing $\mathbb{Q}$ adds a club $C_{\kappa} \subseteq \kappa$ such that $\delta \in C_{\kappa}$ implies $o(\delta)^V < F(\delta)$.  Since $j(F)(\kappa) = F(\kappa)$ and $o(\kappa)^M < F(\kappa)$, the $\kappa$th stage of $j(\mathbb{P}_{\kappa})$ is $\mathbb{Q}$.  Thus, $j(\mathbb{P}_{\kappa})$ factors as $\mathbb{P}_{\kappa} \ast \dot{\mathbb{Q}} \ast \mathbb{P}_{\text{tail}}$ where $\mathbb{P}_{\text{tail}}$ is the forcing past stage $\kappa$.  Since $G_{\kappa} \ast g$ is generic for $\mathbb{P}_{\kappa} \ast \dot{\mathbb{Q}}$, we just need an $M[G_{\kappa}][g]$-generic filter for $\mathbb{P}_{\text{tail}}$ which is in $V[G_{\kappa}][g]$ (since we want the final embedding to be there).  Thus, we shall diagonalize to get an $M[G_{\kappa}][g]$-generic filter for $\mathbb{P}_{\text{tail}}$ after checking we have met the criteria.  Namely, since $|\mathbb{P}_{\kappa}| = |\mathbb{Q}| = \kappa$, the forcing $\mathbb{P}_{\kappa} \ast \dot{\mathbb{Q}}$ has the $\kappa^+$-chain condition.  Since $M^{\kappa} \subseteq M$, it follows [Hamkins1 (Theorem 54)] that $M[G_{\kappa}][g]^{\kappa} \subseteq M[G_{\kappa}][g]$.  Also, since $2^{\kappa} = \kappa^+$, the forcing $\mathbb{P}_{\kappa}$ has $\kappa^+$ many dense sets in $V$.  Thus, $\mathbb{P}_{\text{tail}}$ has at most $|j(\kappa^+)|^V \le \kappa^{+^{\kappa}} = \kappa^+$ many dense sets in $M[G_{\kappa}][g]$.  Also, since for any $\beta < \kappa$, the forcing $\mathbb{P}_{\kappa}$ has a dense subset which is ${\le}\beta$-closed (by Lemma \ref{Lemma. Club shooting.}), it follows that there is a dense subset of $\mathbb{P}_{\text{tail}}$ which is ${\le}\kappa$-closed.  Thus, we may diagonalize to get $G^* \subseteq \mathbb{P}_{\text{tail}}$ for this dense subset, a generic filter in $V[G_{\kappa}][g]$, which is $M[G_{\kappa}][g]$-generic.  Thus, $j(G_{\kappa}) = G_{\kappa} \ast g \ast G^*$ is $M$-generic for $j(\mathbb{P}_{\kappa})$ and $j''G_{\kappa} \subseteq G_{\kappa} \ast g \ast G^*$, so we may lift $j$ to $j:V[G_{\kappa}] \to M[j(G_{\kappa})]$.  

Next, we lift $j$ through $\mathbb{Q}$.  The forcing $j(\mathbb{Q})$ adds a club $C \subseteq j(\kappa)$ such that $\delta \in C$ implies $o(\delta)^{M[j(G_{\kappa})]} < j(F)(\delta)$.  Since $j(G_{\kappa}) \in V[G_{\kappa}][g]$ and $M^{\kappa} \subseteq M$, it follows [Hamkins1 (Theorem 53)] that $M[j(G_{\kappa})]^{\kappa} \subseteq M[j(G_{\kappa})]$.  Also, by Lemma \ref{Lemma. Club shooting.}, there is a dense subset of $j(\mathbb{Q})$ which is ${\le}\kappa$-closed, and since $|\mathbb{Q}| = \kappa$, it follows that $j(\mathbb{Q})$ has at most $\kappa^+$ many dense subsets in $M[j(G_{\kappa})]$.  Let $c = \cup g$ and consider the set $\bar{c} = c \cup \{ \kappa \}$.  By the construction of the forcing, $\forall \delta < \kappa$, $\delta \in c$ implies $o(\delta)^V < F(\delta)$.  Thus, $\forall \delta < \kappa$, $\delta \in c$ implies $o(\delta)^M < F(\delta)$. By Lemma \ref{Lemma. Mitchell rank doesn't go up.} it follows that $\delta \in c$ implies $o(\delta)^{M[j(G)]} < F(\delta)$.  Since $o(\kappa)^{M[j(G)]} \le o(\kappa)^M = \beta_0 < \min \{ o(\kappa)^V, F(\kappa) = \alpha\}$ it follows that $o(\kappa)^{M[j(G)]} < \alpha = j(F)(\kappa)$.  Thus $\bar{c}$ is a closed, bounded subset of $j(\kappa)$ such that if $\delta \in \bar{c}$ then $o(\delta)^{M[j(G)]} < j(F)(\delta)$.  Thus $\bar{c}$ is a condition of $j(\mathbb{Q})$.  Thus, diagonalize to get an $M[j(G)]$-generic filter, $g^* \subseteq j(\mathbb{Q})$ which contains the condition $\bar{c}$, so that $j''g \subseteq j(g) = g^*$.  Therefore, we may lift the embedding to $j: V[G_{\kappa}][g] \to M[j(G_{\kappa})][j(g)]$.

From the previous lifting arguments will follow that $o(\kappa)^{V[G]} \ge \min \{ o(\kappa)^V, F(\kappa) \}$.  The lifted embedding $j:V \to M$ was chosen so that $o(\kappa)^M = \beta_0 < \min \{ o(\kappa)^V, F(\kappa) \}$. Since $\beta_0$ was arbitrary, by Lemma \ref{Lemma. Mitchell rank doesn't go up.} it follows that $o(\kappa)^{V[G]} \ge \min \{ o(\kappa)^V, F(\kappa) \}$.  All that is left is to see that $o(\kappa)^{V[G]} \le \min \{o(\kappa)^V, F(\kappa) \}$.  By Lemma \ref{Lemma. Mitchell rank doesn't go up.} which shows that Mitchell rank does not go up for $V \subseteq V[G]$, we have $o(\kappa)^{V[G]} \le o(\kappa)^V$.  Consider any embedding $j: V[G] \to M[j(G)]$ in $V[G]$, and the new club $c \subseteq \kappa$ where $\forall \delta < \kappa$, $\delta \in c$ implies $o(\delta)^V < F(\delta)$.  By the induction hypothesis, $\delta \in c$ implies $o(\delta)^{V[G]} < F(\delta)$.  Apply $j$ to this statement to get $\forall \delta < j(\kappa)$, $\delta \in j(c)$ implies $o(\delta)^{M[j(G)]} < j(F)(\kappa) = F(\kappa)$.  Note that since $F \restriction \kappa$ represents $F(\kappa)$ in $V[G]$ and $M[j(G)]^{\kappa} \subseteq M[j(G)]$, it follows that $F \restriction \kappa$ represents $F(\kappa)$ in $M[j(G)]$ as well.  Since $j$ was arbitrary, $\kappa < j(\kappa)$, and $\kappa \in j(c)$ it now follows by Lemma \ref{Lemma. Mitchell rank in M.} it follows that $o(\kappa)^{V[G]} \le F(\kappa)$.  Thus $o(\kappa)^{V[G]} \le \min \{ o(\kappa)^V, F(\kappa) \}$.  Thus $o(\kappa)^{V[G]} = \min \{ o(\kappa)^V, F(\kappa) \}$.  

\end{proof}
\end{theorem}

\newpage

\section{\textbf{Supercompact and strongly compact cardinals.} }

This section begins with supercompact and strongly compact cardinals.  The main theorem is about forcing to change degrees of supercompact cardinals.  An uncountable cardinal $\kappa$ is $\theta$-\textit{supercompact} if and only if there is an elementary embedding $j:V \to M$, with critical point $\kappa$ and $M^{\theta} \subseteq M$, where $\kappa < \theta < j(\kappa)$.  Equivalently, $\kappa$ is $\theta$-supercompact if and only if there is a normal fine measure on $P_{\kappa} \theta$.   A \textit{fine measure} in the previous sentence means a $\kappa$-complete ultrafilter such that $\forall \sigma \in P_{\kappa}\theta$ the set $\{ X \subseteq P_{\kappa} \theta \ | \ \sigma \in X \}$ is in the filter.  A cardinal $\kappa$ is \textit{supercompact} if and only if it is $\theta$-supercompact for every $\theta > \kappa$.  Notice that $\kappa$ measurable means that $\kappa$ is $\kappa$-supercompact.  An uncountable cardinal $\kappa$ is $\theta$-\textit{strongly compact} if and only if there is an elementary embedding $j: V \to M$, with critical point $\kappa$, such that for all $t \subseteq M$, where $|t|^V \le \theta$, there is an $s \in M$ such that $t \subseteq s$ and $|s|^M < j(\kappa)$.  Equivalently, a cardinal $\kappa$ is $\theta$-strongly compact for $\kappa < \theta$ if and only if there is a $\kappa$-complete fine measure on $P_{\kappa} \theta$.  An uncountable cardinal $\kappa$ is \textit{strongly compact} if and only if it is $\theta$-strongly compact, for every $\theta > \kappa$.  Equivalently, an uncountable regular cardinal $\kappa$ is strongly compact if and only if for any set $S$, every $\kappa$-complete filter on $S$ can be extended to a $\kappa$-complete ultrafilter on $S$ [Jech p. 365].  The following theorems show how to distinguish, by forcing, between levels of supercompactness, levels of strong compactness, between measurable and supercompact cardinals, between measurable and strongly compact cardinals, and between strongly compact and supercompact cardinals.  

\begin{theorem} \label{Theorem. Kill supercompact softly.} If $\kappa$ is $<\theta$-supercompact, for $\theta \ge \kappa$ regular with $\beta^{<\kappa} < \theta$ for all $\beta < \theta$, then there is a forcing extension where $\kappa$ is $<\theta$-supercompact, but not $\theta$-supercompact, and indeed not even $\theta$-strongly compact.

\begin{proof} Suppose $\kappa$ is $<\theta$-supercompact for some regular $\theta > \kappa$, and let $\mathbb{P} = Add(\omega, 1) \ast Add(\theta, 1)$.  By [Hamkins/Shelah], after forcing with $\mathbb{P}$, the cardinal $\kappa$ is not $\theta$-supercompact.  However, adding a Cohen real is small relative to $\kappa$, so the first factor of $\mathbb{P}$ preserves that $\kappa$ is $<\theta$-supercompact.  The second factor of $\mathbb{P}$ is $<\theta$-closed, so that a normal fine measure on $P_{\kappa} \beta$, where $\beta < \theta$ is still a normal fine measure on $P_{\kappa} \beta$, since no new subsets of $\beta$ were added.  Let $G \subseteq \mathbb{P}$ be $V$-generic. Thus, in $V[G]$, the cardinal $\kappa$ is $<\theta$-supercompact, but not $\theta$-supercompact.  

\end{proof}
\end{theorem}

\begin{corollary} If $\kappa$ is $\theta$-supercompact where $\kappa < \theta$ and $\theta^{<\kappa} = \theta$, then there is a forcing extension where $\kappa$ is $\theta$-supercompact, but not $\theta^+$-supercompact.  

\end{corollary}

\begin{theorem} If $\kappa$ is ${<}\theta$-strongly compact, for $\theta \ge \kappa$ regular with $\beta^{<\kappa} < \theta$ for all $\beta < \theta$, then there is a forcing extension where $\kappa$ is ${<}\theta$-strongly compact, but not $\theta$-strongly compact.  

\begin{proof} Suppose $\kappa$ is $<\theta$-strongly compact for regular $\theta \ge \kappa$, and let $\mathbb{P} = Add(\omega, 1) \ast Add(\theta, 1)$.  By [Hamkins/Shelah], this forcing destroys the $\theta$-strong compactness of $\kappa$.  However, adding a Cohen real is small relative to $\kappa$, so the first factor preserves that $\kappa$ is $<\theta$-strongly compact.  Also, since the second factor is $<\theta$-closed, any $\kappa$-complete fine measure on $P_{\kappa}\beta$, where $\beta < \theta$, is still a $\kappa$-complete fine measure since there are no new subsets of $\beta$ to measure.  Thus, $\mathbb{P}$ preserves that $\kappa$ is $<\theta$-strongly compact.  Let $G \subseteq \mathbb{P}$ be $V$-generic. Then, in $V[G]$, the cardinal $\kappa$ is $<\theta$-strongly compact, but not $\theta$-strongly compact.  
\end{proof}
\end{theorem}

\begin{corollary}  If $\kappa$ is $\theta$-strongly compact, where $\kappa \le \theta$ and $\theta^{<\kappa} = \theta$, then there is a forcing extension where $\kappa$ is $\theta$-strongly compact, but not $\theta^+$-strongly compact.  

\end{corollary}

\begin{corollary} If $\kappa$ is measurable, then there is a forcing extension where $\kappa$ is measurable, but not strongly compact.

\begin{proof} Let $\theta = \kappa$  in the proof of Theorem \ref{Theorem. Kill supercompact softly.}
\end{proof}
\end{corollary}

The following theorem is due to Magidor [Magidor]. It shows how to force between strong compactness and supercompactness.  

\begin{theorem} (Magidor) If $\kappa$ is strongly compact, then there is a forcing extension where $\kappa$ is strongly compact, but not supercompact.  
\end{theorem}

A word on the proof of the theorem.  Suppose $\kappa$ is a strongly compact cardinal.  Then, there is a forcing extension, $V[G]$, where $\kappa$ is the least strongly compact cardinal, and the least measurable cardinal [Magidor].  However, the least measurable cardinal can never be supercompact, since there are many measurable cardinals below any supercompact cardinal.  Thus, in $V[G]$, the cardinal $\kappa$ is strongly compact, but not supercompact.

As for measurable cardinals, an analgous rank to Mitchell rank can be assigned to supercompactness embeddings.  Suppose a cardinal $\kappa$ is $\theta$-supercompact for fixed $\theta$.  If $\mu$ and $\nu$ are normal fine measures on $P_{\kappa}{\theta}$, define the Mitchell relation $\mu \lhd_{\theta \text{-sc}} \nu$ if and only if $\mu \in M_{\nu}$ where $j:V \to M_{\nu}$ is an ultrapower embedding by $\nu$.  Since the relation $\lhd_{\theta \text{-sc}}$ is well-founded for given $\kappa$ and $\theta$, the Mitchell rank of $\mu$ is its rank with respect to $\lhd_{\theta \text{-sc}}$. Thus, for fixed $\kappa$ and $\theta$, the notation $o_{sc}(\kappa) = \alpha$ means the height of $\lhd_{\theta \text{-sc}}$ on normal fine measures on $P_{\kappa}\theta$ is $\alpha$. By definition, $o_{\theta \text{-sc}}(\kappa)$ is the height of the well-founded Mitchell relation $\lhd_{\theta \text{-sc}}$ on $m(\kappa)$, the collection of normal fine measures on $P_{\kappa} \theta$. Thus the definition of rank for $\mu \in m(\kappa)$ with respect to $\lhd_{\theta \text{-sc}}$ is $o_{\theta \text{-sc}}(\mu) = \sup \{ o_{\theta \text{-sc}}(\nu) + 1 \ | \nu \lhd_{\theta \text{-sc}} \mu \}$. Thus $o_{\theta \text{-sc}}(\kappa) = \{ o_{\theta \text{-sc}}(\mu) + 1 \ | \ \mu \in m(\kappa) \}$. So that $o_{\theta \text{-sc}}(\kappa) = 0$ means that $\kappa$ is not $\theta$-supercompact, and $o_{\theta \text{-sc}}(\kappa) > 0$ means that $\kappa$ is $\theta$-supercompact.  

The main theorems of this section will rely on the following Lemmas about Mitchell rank for supercompactness and are analogues of Lemmas \ref{Lemma. Mitchell rank in M.}, \ref{Lemma. Mitchell rank doesn't go up.}, and \ref{Lemma. Lift normal ultrapower embeddings.} for Mitchell rank.

\begin{lemma} \label{Lemma. Mitchell rank for supercompactness in M.}  

\begin{proof}
Suppose $\kappa$ is $\theta$-supercompact, $\mu$ is a normal fine measure on $P_{\kappa} \theta$, and $j_{\mu}: V \to M_{\mu}$ is the ultrapower emebedding by $\mu$ with critical point $\kappa$. The model $M_{\mu}$ computes $\lhd_{\theta \text{-sc}}$ correctly for $\nu \lhd_{\theta \text{-sc}} \mu$, which means $\nu \in M_{\mu}$, since $M^{\theta} \subseteq M$, and so $M_{\mu}$ has  all the functions needed to check whether $\nu \lhd_{\theta \text{-sc}} \mu$.  Then
\begin{align*}
    o_{\theta \text{-sc}}(\kappa)^{M_{\mu}}  &= \sup\{ o_{\theta \text{-sc}}(\nu)^{M_{\nu}} + 1 \ | \ \nu \in M_{\mu} \} \\ 
     &= \sup\{o_{\theta \text{-sc}}(\nu)^{M_{\mu}} + 1 \ | \ \nu \lhd_{\theta \text{-sc}} \mu \} \\
     &= \sup\{o_{\theta \text{-sc}}(\nu) + 1 \ | \ \nu \lhd_{\theta \text{-sc}} \mu \} \\
     &= o_{\theta \text{-sc}}(\mu)
\end{align*}
Thus, the Mitchell rank for supercompactness of $\kappa$ in $M_{\mu}$ is the Mitchell rank of $\mu$.  Let $\beta < o_{\theta \text{-sc}}(\kappa)$.  Since $\lhd_{\theta \text{-sc}}$ is well-founded, there is $\mu$ a normal fine measure on $P_{\kappa} \theta$ such that $o_{\theta \text{-sc}}(\mu) = \beta$.  Then if $j_{\mu} : V \to M_{\mu}$ is the ultrapower by $\mu$ with critical point $\kappa$, then $o_{\theta \text{-sc}}(\kappa)^{M_{\mu}} = o_{\theta \text{-sc}}(\mu) = \beta$.  Conversely, suppose there is an elementary embedding $j:V \to M$ with cp$(j) = \kappa$ and $o_{\theta \text{-sc}}(\kappa)^M = \beta$.  Let $\mu$ be the induced normal fine measure, using $j''\theta$ as a seed ($X \in \mu$ if and only if $j''\theta \in j(X)$) and let $M_{\mu}$ be the ultrapower by $\mu$.  Then $M$ and $M_{\mu}$ have the same normal fine measures and the same Mitchell order and ranks for $\theta$-supercompactness, so $M_{\mu} \models o_{\theta \text{-sc}}(\kappa) = \beta$.  Since we have established $o_{\theta \text{-sc}}(\kappa)^{M_{\mu}} = o_{\theta \text{-sc}}(\mu)$, the Mitchell rank for supercompactness of $\mu$ in $m(\kappa)$ is $\beta$, that is $o_{\theta \text{-sc}}(\mu) = \beta$.  Thus since $o_{\theta \text{-sc}}(\kappa) = \{ o_{\theta \text{-sc}}(\mu) + 1 \ | \ \mu \in m(\kappa) \}$, it follows that $o_{\theta \text{-sc}}(\kappa)^V > \beta$.

\end{proof}
\end{lemma}
The following Lemma is an analogue of Lemma \ref{Lemma. Mitchell rank doesn't go up.} for degrees of Mitchell rank for supercompactness.  It is a generalization of Corollary 26 in [Hamkins4] which states that if $V \subseteq \bar{V}$ satisfies the approximation and cover properties, then no supercompact cardinals are created from $V$ to $\bar{V}$.  

\begin{lemma} \label{Lemma. Mitchell rank for supercompactness doesn't go up.} Suppose $V \subseteq \bar{V}$ satisfies the $\delta$ approximation and cover properties.  If $\kappa, \theta > \delta$, then $o_{\theta \text{-sc}}(\kappa)^{\bar{V}} \le o_{\theta \text{-sc}}(\kappa)^{V}$.

\begin{proof}
Suppose $V \subseteq \bar{V}$ satisfies the $\delta$ approximation and cover properties.  Suppose $\kappa, \theta > \delta$.  We want to show $o_{\theta \text{-sc}}(\kappa)^{\bar{V}} \le o_{\theta \text{-sc}}(\kappa)^V$ by induction on $\kappa$.  Assume inductively that for all $\kappa' < \kappa$, we have $o_{\theta \text{-sc}}(\kappa')^{\bar{V}} \le o_{\theta \text{-sc}}(\kappa')^V$.  If $o_{\theta \text{-sc}}(\kappa)^{\bar{V}} > o_{\theta \text{-sc}}(\kappa)^V$ then by Lemma \ref{Lemma. Mitchell rank for supercompactness in M.} get a $\theta$-supercompactness embedding $j: \bar{V} \to \bar{M}$ with critical point $\kappa$ and $o_{\theta \text{-sc}}(\kappa)^{\bar{M}} = o_{\theta \text{-sc}}(\kappa)^{V}$.   By [Hamkins4], $j$ is a lift of a $\theta$-supercompactness embedding $j \restriction V: V \to M$ for some $M$, and this embedding is a class in $V$.  Again by \ref{Lemma. Mitchell rank for supercompactness in M.} we have established that $o_{\theta \text{-sc}}(\kappa)^M < o_{\theta \text{-sc}}(\kappa)^V$.  Thus, $o_{\theta \text{-sc}}(\kappa)^M < o_{\theta \text{-sc}}(\kappa)^{\bar{M}}$.  However, this contradicts that $j$ applied to the induction hypothesis gives $o_{\theta \text{-sc}}(\kappa)^{\bar{M}} \le o_{\theta \text{-sc}}(\kappa)^{M}$ since $\kappa < j(\kappa)$.  Therefore, $o_{\theta \text{-sc}}(\kappa)^{\bar{V}} \le o_{\theta \text{-sc}}(\kappa)^V$.
\end{proof}
\end{lemma}

\begin{lemma} \label{Lemma. Lift normal ultrapower embeddings for supercompactness.} If $V \subseteq V[G]$ is a forcing extension and every $\theta$-supercompactness embedding $j:V \to M$ in $V$ (with any critical point $\kappa$) lifts to an embedding $j:V[G] \to M[j(G)]$, then $o_{\theta \text{-sc}}(\kappa)^{V[G]} \ge o_{\theta \text{-sc}}(\kappa)^V$.  In addition, if $V \subseteq V[G]$ also has $\delta$ approximation and $\delta$ cover properties, then $o_{\theta \text{-sc}}(\kappa)^{V[G]} = o_{\theta \text{-sc}}(\kappa)^V$.
\begin{proof} Suppose $V \subseteq V[G]$ is a forcing extension where every normal ultrapower embedding in $V$ lifts to a normal ultrapower embedding in $V[G]$.  Suppose $\kappa$ is a $\theta$-supercompact cardinal, and for every $\theta$-supercompact cardinal $\kappa' < \kappa$, assume $o_{\theta \text{-sc}}(\kappa')^{V[G]} \ge o_{\theta \text{-sc}}(\kappa')^V$.  Fix any $\beta < o_{\theta \text{-sc}}(\kappa)$.  By Lemma \ref{Lemma. Mitchell rank for supercompactness in M.} there exists an elementary embedding $j:V \to M$ with $cp(j) = \kappa$ and $M \models o_{\theta \text{-sc}}(\kappa) = \beta$.  Then, by the assumption on $V \subseteq V[G]$, this $j$ lifts to $j:V[G] \to M[j(G)]$.  Then, applying $j$ to the induction hypothesis gives $o_{\theta \text{-sc}}(\kappa)^{M[j(G)]} \ge o_{\theta \text{-sc}}(\kappa)^M$ since $\kappa < j(\kappa)$.  Then, since $M \models o_{\theta \text{-sc}}(\kappa) = \beta$ it follows that $o_{\theta \text{-sc}}(\kappa)^{M[j(G)]} \ge \beta$.  Thus by Lemma \ref{Lemma. Mitchell rank for supercompactness in M.} we have $o_{\theta \text{-sc}}(\kappa)^{V[G]} > \beta$.  Thus $o_{\theta \text{-sc}}(\kappa)^{V[G]} \ge o_{\theta \text{-sc}}(\kappa)^V$.  

If $V \subseteq V[G]$ also satisfies the $\delta$ approximation and $\delta$ cover properties, then by Lemma \ref{Lemma. Mitchell rank for supercompactness doesn't go up.} the Mitchell rank does not go up between these models. Thus $_{\theta \text{-sc}}o(\kappa)^{V[G]} = o_{\theta \text{-sc}}(\kappa)^V$.
\end{proof}

\end{lemma}
The following theorem shows how to force a $\kappa^+$-supercompact cardinal $\kappa$ to have Mitchell rank for supercompactness at most 1.

\begin{theorem}  If $V \models ZFC + GCH$ then there is a forcing extension where every cardinal $\kappa$ which is $\kappa^+$-supercompact has $o_{\kappa^+ \text{-sc}}(\kappa)^{V[G]} = \min \{o_{\kappa^+ \text{-sc}}(\kappa)^V, 1 \}$.  

\begin{proof}
Suppose $V \models ZFC + GCH$. Let $\mathbb{P}$ be an Easton support Ord-length iteration which forces at inaccessible stages $\gamma$ to add a club $c_{\gamma} \subseteq \gamma$ such that $\delta \in c_{\gamma}$ implies $o_{\delta^+ \text{-sc}}(\delta)^V = 0$.  Let $G \subseteq \mathbb{P}$ be $V$-generic.  Let $\delta_0$ be the first inaccessible cardinal.  The forcing before stage $\delta_0$ is trivial, and the forcing at stage $\delta_0$ adds a club to $\delta_0$.  Thus the forcing $\mathbb{P}$ up to and including stage $\delta_0$ has cardinality $\delta_0$, and the forcing past stage $\delta_0$ has a dense set which is closed up to the next inaccessible cardinal.  Thus $\mathbb{P}$ has a closure point at $\delta_0$ which implies $V \subseteq V[G]$ satisfies the $\delta_0^+$ approximation and cover properties by [Hamkins4 (Lemma 13)].  Thus by [Hamkins4 (Corollary 26)], the forcing $\mathbb{P}$ does not create $\theta$-supercompact cardinals for any $\theta$.  Note that any $\theta$-supercompact cardinal $\kappa$, both $\theta$ and $\kappa$ are above this closure point.  

Suppose $\kappa$ is $\kappa^+$-supercompact.  The induction hypothesis is that any $\kappa' < \kappa$ has $o_{\kappa'^+ \text{-sc}}(\kappa')^{V[G]} = \min \{o(\kappa')^V, 1\}$.  Pick $j:V \to M$ a $\kappa^+$-supercompactness embedding with critical point $\kappa$ such that $o_{\kappa^+ \text{-sc}}(\kappa)^M = 0 < \min\{ o_{\kappa^+ \text{-sc}}(\kappa)^V, 1\}$.  We shall lift $j$ through $\mathbb{P}$.  The forcing above $\kappa$ does not affect the Mitchell rank for $\kappa^+$-supercompactness of $\kappa$ since the forcing past stage $\kappa$ has a dense subset which is closed up to the next inaccessible cardinal past $\kappa$.  Call $\mathbb{P}_{\kappa} \ast \dot{\mathbb{Q}}$ the forcing up to and including stage $\kappa$.  First we shall lift $j$ through $\mathbb{P}_{\kappa}$.  Let $G_{\kappa} \subseteq \mathbb{P}_{\kappa}$ be $V$-generic and $g \subseteq \mathbb{Q}$ be $V[G_{\kappa}]$-generic.  First, we need an $M$-generic filter for $j(\mathbb{P}_{\kappa})$.  Since $cp(j)= \kappa$, the forcings $\mathbb{P}_{\kappa}$ and $j(\mathbb{P}_{\kappa})$ agree up to stage $\kappa$.  The forcing $\mathbb{Q}$ adds a club $c \subseteq \kappa$ such that $\forall \delta < \kappa$, $\delta \in c$ implies $o_{\delta^+ \text{-sc}}(\delta)^V = 0$.  Since $o_{\delta^+ \text{-sc}}(\delta)^V = 0$ implies $o_{\delta^+ \text{-sc}}(\delta)^M = 0$, it follows that the $\kappa$th stage of $j(\mathbb{P}_{\kappa})$ is $\mathbb{Q}$.  Thus, $j(\mathbb{P}_{\kappa})$ factors as $\mathbb{P}_{\kappa} \ast \dot{\mathbb{Q}} \ast \mathbb{P}_{\text{tail}}$ where $\mathbb{P}_{\text{tail}}$ is the forcing past stage $\kappa$.  Since $G_{\kappa} \ast g$ is $M$-generic for $\mathbb{P}_{\kappa} \ast \dot{\mathbb{Q}}$, we only need an $M[G_{\kappa}][g]$-generic filter for $\mathbb{P}_{\text{tail}}$ which is $V[G_{\kappa}][g]$.  We will diagonalize to obtain such a filter, after checking that we have met the criteria.  Since $|\mathbb{P}_{\kappa}| = |\mathbb{Q}| = \kappa$, and $V \models \text{GCH}$, it follows that the forcing $\mathbb{P}_{\kappa} \ast \dot{\mathbb{Q}}$ has the $\kappa^+$-chain condition.  Since $M^{\kappa^+} \subseteq M$, it follows [Hamkins1 (Theorem 54)] that $M[G_{\kappa}][g]^{\kappa^+} \subseteq M[G_{\kappa}][g]$.  Since $|\mathbb{P}_{\kappa}| = \kappa$, and has the $\kappa^+$-chain condition, it has at most $\kappa^{<\kappa^+} = \kappa^+$ many maximal antichains.  Since $|j(\kappa^+)| \le \kappa^{+^{\kappa^{+^{< \kappa}}}} = \kappa^{++}$, it follows that $\mathbb{P}_{\text{tail}}$ has at most $\kappa^{++}$ many maximal antichains in $M[G_{\kappa}][g]$.  Since $\mathbb{P}_{\text{tail}}$ has a dense subset which is closed up to the next inaccessible past $\kappa$, it has a dense subset which is ${\le}\kappa^+$-closed.  Thus, we can enumerate the maximal antichains of $\mathbb{P}_{\text{tail}}$, and build a descending sequence of conditions meeting them all through the dense set which has ${\le}\kappa^+$-closure, using this closure and the fact that $M[G_{\kappa}][g]$ is closed under $\kappa^+$-sequences to build through the limit stages.  Then the upward closure of this sequence $G_{\text{tail}} \subseteq \mathbb{P}_{\text{tail}}$ is $M[G_{\kappa}][g]$-generic and is in $V[G_{\kappa}][g]$.  Thus $j(G_{\kappa}) = G_{\kappa} \ast g \ast G_{\text{tail}}$ is $M$-generic, it is in $V[G_{\kappa}][g]$ and $j''G_{\kappa} \subseteq j(G_{\kappa})$.  Thus we have the partial lift $j:V[G_{\kappa}] \to M[j(G_{\kappa})]$.

Next we lift through $\mathbb{Q}$ by finding an $M[j(G_{\kappa})]$-generic filter for $j(\mathbb{Q})$.  The forcing $j(\mathbb{Q})$ adds a club $C \subseteq j(\kappa)$ such that $\forall \delta < j(\kappa)$, $\delta \in C$ implies $o_{\delta^+ \text{-sc}}(\delta)^{M[j(G_{\kappa})]} = 0$.  Since $j(G_{\kappa}) \in V[G_{\kappa}][g]$ and $M^{\kappa^+} \subseteq M$, it follows that $M[j(G_{\kappa})]^{\kappa^+} \subseteq M[j(G_{\kappa})]$ by [Hamkins1 (Theorem 53)].  Since for every $\beta < \kappa$, there is a dense subset of $\mathbb{Q}$ which is ${\le}\beta$-closed, it follows that there is a dense subset of $j(\mathbb{Q})$ which is ${\le}\kappa^+$-closed.  Also since $|\mathbb{Q}| = \kappa$ and $\mathbb{Q}$ has the $\kappa^+$-chain condition, the forcing $\mathbb{Q}$ has at most $\kappa^+$ many maximal antichains.  Therefore $j(\mathbb{Q})$ has at most $|j(\kappa^+)| \le \kappa^{+^{\kappa^{+^{< \kappa}}}} = \kappa^{++}$ many maximal antichains in $M[j(G_{\kappa})]$.  Thus we can diagonalize to get a generic filter, but we will do so with a master condition as follows.  Let $c = \cup g$ the new club of $\kappa$ and consider $\bar{c} = c \cup \{ \kappa \}$ which is in $M[j(G_{\kappa})]$.  We have $\delta \in c$ implies $o_{\delta^+ \text{-sc}}(\delta)^V = 0$.  Thus, $\delta \in c$ implies $o_{\delta^+ \text{-sc}}(\delta)^M = 0$.  By Corollary 26 [Hamkins4], no new $\delta^+$-supercompact cardinals are created between $M$ and $M[j(G_{\kappa})]$, so that $\delta \in c$ implies $o_{\delta^+ \text{-sc}}(\delta)^{M[j(G_{\kappa})]} = 0$.  Since $o_{\kappa^+ \text{-sc}}(\kappa)^{M} = 0$, it follows that $o_{\kappa^+ \text{-sc}}(\kappa)^{M[j(G_{\kappa})]} = 0$.  Thus $\bar{c}$ is a closed, bounded subset of $j(\kappa)$ such that every $\delta$ in $\bar{c}$ has $o_{\delta^+ \text{-sc}}(\delta)^{M[j(G_{\kappa})]} = 0$.  Thus $\bar{c} \in j(\mathbb{Q})$.  Thus diagonalize to get an $M[j(G_{\kappa})]$-generic filter $g^* \subseteq j(\mathbb{Q})$ which contains the master condition $\bar{c}$.  Then $j''g \subseteq j(g) = g^*$.  Therefore we may lift the embedding to $j:V[G_{\kappa}][g] \to M[j(G_{\kappa})][j(g)]$.  

The success of the lifting arguments show that $o_{\kappa^+ \text{-sc}}(\kappa)^{V[G]} \ge \min\{o_{\kappa^+ \text{-sc}}(\kappa)^V, 1\}$.  By Lemma \ref{Lemma. Mitchell rank for supercompactness doesn't go up.} we have $o_{\kappa^+ \text{-sc}}(\kappa)^{V[G]} \le o_{\kappa^+ \text{-sc}}(\kappa)^V$.  All that remains is to see that $o_{\kappa^+ \text{-sc}}(\kappa)^{V[G]} \le 1$. Consider any $\kappa^+$-supercompact embedding $j:V[G] \to M[j(G)]$ with critical point $\kappa$ in $V[G]$, and the new club $c \subseteq \kappa$, where $\forall \delta < \kappa$, $\delta \in c$ implies $o_{\delta^+ \text{-sc}}(\delta)^V = 0$.  This club is in any normal fine measure on $P_{\kappa} \kappa^+$.  By the induction hypothesis, $\forall \delta < \kappa$, $\delta \in c$ implies $o_{\delta^+ \text{-sc}}(\delta)^{V[G]} = 0$.  Applying $j$ to this statement gives $\forall \delta < j(\kappa)$, $\delta \in j(c)$ implies $o_{\delta^+ \text{-sc}}(\delta)^{M[j(G)]} = 0$.  Thus by Lemma \ref{Lemma. Mitchell rank for supercompactness in M.}, since $j$ was arbitrary, $\kappa < j(\kappa)$, and $\kappa \in j(c)$ it follows that $o_{\kappa^+ \text{-sc}}(\kappa)^{M[j(G)]} = 0$, and hence $o_{\kappa^+ \text{-sc}}(\kappa)^{V[G]} \le 1$.  Thus $o_{\kappa^+ \text{-sc}}(\kappa)^{V[G]} \le \min \{ o_{\kappa^+ \text{-sc}}(\kappa)^V, 1 \}$.  Thus for any $\kappa$ which is $\kappa^+$-supercompact, $o_{\kappa^+ \text{-sc}}(\kappa)^{V[G]} = \min \{ o_{\kappa^+ \text{-sc}}(\kappa)^V, 1 \}$.  
\end{proof}
\end{theorem}

For the following theorem for changing Mitchell rank for supercompactness, we will again need representing functions because I would like to now consider cases where $o_{\kappa^+ \text{-sc}}(\kappa) \ge \kappa$. Then for the most general theorem I would like to consider the cases where $\kappa$ is $\theta$-supercompact for $\theta > \kappa^+$. In [Hamkins2] representing functions are considered for supercompactness embeddings as well. Let $\alpha \in H_{\theta^+}$ and $\kappa$ be a $\theta$-supercompact cardinal. The ordinal $\alpha$ is represented by $f: \kappa \to V_{\kappa}$ if whenever $j: V \to M$ an elementary embedding with cp$(j) = \kappa$ and $M^{\theta} \subseteq M$, then $j(f)(\kappa) = \alpha$ (this can be formalized as a first-order statement using extenders) [Hamkins2].  Let $F: \text{Ord} \to \text{Ord}$. The next theorem shows how to change the Mitchell rank for supercompactness for all $\gamma^+$-supercompact cardinals $\gamma$ for which $F \restriction \gamma$ represents $F(\gamma)$.

\begin{theorem} For any $V \models ZFC + GCH$ and any $F: \text{Ord} \to \text{Ord}$, there is a forcing extension $V[G]$ where, if $\kappa$ is $\kappa^+$-supercompact and $F \restriction \kappa$ represents $F(\kappa)$ in V, then $o_{\kappa^+ \text{-sc}}(\kappa)^{V[G]} =\min\{o_{\kappa^+ \text{-sc}}(\kappa)^V, F(\kappa)\}$.

\begin{proof}
Suppose $V \models ZFC + GCH$. Let $F : \text{Ord} \to \text{Ord}$.  Let $\mathbb{P}$ be an Easton support Ord-length iteration which forces at inaccessible stages $\gamma$ to add a club $c_{\gamma} \subseteq \gamma$ such that $\forall \delta < \gamma$, $\delta \in c_{\gamma}$ implies $o_{\delta^+ \text{-sc}}(\delta)^V < F(\delta)$.  Let $G \subseteq \mathbb{P}$ be $V$-generic.  Let $\delta_0$ be the first inaccessible cardinal.  The forcing before stage $\delta_0$ is trivial, and the forcing at stage $\delta_0$ adds a club to $\delta_0$.  Thus the forcing $\mathbb{P}$ up to and including stage $\delta_0$ has cardinality $\delta_0$ and the forcing past stage $\delta_0$ has a dense set which is closed up to the next inaccessible cardinal.  Thus $\mathbb{P}$ has a closure point at $\delta_0$, thus $V \subseteq V[G]$ satisfies the $\delta_0^+$ approximation and cover properties by [Hamkins4 (Lemma 13)].  Thus by [Hamkins4 (Corollary 26)], the forcing $\mathbb{P}$ does not create $\theta$-supercompact cardinals for any $\theta$.  Note that for any $\theta$-supercompact cardinal $\kappa$, both $\kappa$ and $\theta$ are above this closure point.  

Let $\kappa$ be a $\kappa^+$-supercompact cardinal such that $F \restriction \kappa$ represents $F(\kappa) = \alpha$ in $V$. Here we will need that representing functions are still representing functions in the extension.  The fact we will use for this is from [Hamkins4] which says that if $V \subseteq V[G]$ satisfies the hypothesis of approximation and cover theorem, then every $\theta$-closed embedding $j:V[G] \to M[j(G)]$ is a lift of an embedding $j \restriction V : V \to M$ that is in $V$ (and $j \restriction V$ is also a $\theta$-closed embedding).  Thus, if $F \restriction \kappa$ is a representing function for $F(\kappa)$ in the ground model, $(j \restriction V) (F \restriction \kappa)(\kappa) = F(\kappa)$, and thus $F \restriction \kappa$ is a representing function for $F(\kappa)$ with respect to such embeddings $j$ in $V[G]$.

The induction hypothesis is that any $\kappa' < \kappa$, where $F \restriction \kappa'$ represents $F(\kappa')$ in $V$, has $o_{\kappa'^+ \text{-sc}}(\kappa')^{V[G]} = \min \{o_{\kappa'^+ \text{-sc}}(\kappa')^V, F(\kappa')\}$.  Pick $j:V \to M$ a $\kappa^+$-supercompactness embedding with critical point $\kappa$ such that $o_{\kappa^+ \text{-sc}}(\kappa)^M = \beta_0 < \min\{ o_{\kappa^+ \text{-sc}}(\kappa)^V, F(\kappa)\}$.  We shall lift $j$ through $\mathbb{P}$.  The forcing above $\kappa$ does not affect the Mitchell rank for $\kappa^+$-supercompactness of $\kappa$ since the forcing past stage $\kappa$ has a dense subset which is closed up to the next inaccessible cardinal past $\kappa$.  Call $\mathbb{P}_{\kappa} \ast \dot{\mathbb{Q}}$ the forcing up to and including stage $\kappa$.  First we shall lift $j$ through $\mathbb{P}_{\kappa}$.  Let $G_{\kappa} \subseteq \mathbb{P}_{\kappa}$ be $V$-generic and $g \subseteq \mathbb{Q}$ be $V[G_{\kappa}]$-generic.  First, we need an $M$-generic filter for $j(\mathbb{P}_{\kappa})$.  Since $cp(j)= \kappa$, the forcings $\mathbb{P}_{\kappa}$ and $j(\mathbb{P}_{\kappa})$ agree up to stage $\kappa$.  The forcing $\mathbb{Q}$ adds a club $c \subseteq \kappa$ such that $\forall \delta < \kappa$, $\delta \in c$ implies $o_{\delta^+ \text{-sc}}(\delta)^V < F(\delta)$.  Since $o_{\delta^+ \text{-sc}}(\delta)^V < F(\delta)$ implies $o_{\delta^+ \text{-sc}}(\delta)^M < F(\delta) = j(F)(\delta)$, it follows that the $\kappa$th stage of $j(\mathbb{P}_{\kappa})$ is $\mathbb{Q}$.  Thus, $j(\mathbb{P}_{\kappa})$ factors as $\mathbb{P}_{\kappa} \ast \dot{\mathbb{Q}} \ast \mathbb{P}_{\text{tail}}$ where $\mathbb{P}_{\text{tail}}$ is the forcing past stage $\kappa$.  Since $G_{\kappa} \ast g$ is $M$-generic for $\mathbb{P}_{\kappa} \ast \dot{\mathbb{Q}}$, we only need an $M[G_{\kappa}][g]$-generic filter for $\mathbb{P}_{\text{tail}}$ which is $V[G_{\kappa}][g]$.  We will diagonalize to obtain such a filter, after checking that we have met the criteria.  Since $|\mathbb{P}_{\kappa}| = |\mathbb{Q}| = \kappa$, and $V \models \text{GCH}$, it follows that the forcing $\mathbb{P}_{\kappa} \ast \dot{\mathbb{Q}}$ has the $\kappa^+$-chain condition.  Since $M^{\kappa^+} \subseteq M$, it follows [Hamkins1 (Theorem 54)] that $M[G_{\kappa}][g]^{\kappa^+} \subseteq M[G_{\kappa}][g]$.  Since $|\mathbb{P}_{\kappa}| = \kappa$, and has the $\kappa^+$-chain condition, it has at most $\kappa^{<\kappa^+} = \kappa^+$ many maximal antichains.  Since $|j(\kappa^+)| \le \kappa^{+^{\kappa^{+^{< \kappa}}}} = \kappa^{++}$, it follows that $\mathbb{P}_{\text{tail}}$ has at most $\kappa^{++}$ many maximal antichains in $M[G_{\kappa}][g]$.  Since $\mathbb{P}_{\text{tail}}$ has a dense subset which is closed up to the next inaccessible past $\kappa$, it has a dense subset which is ${\le}\kappa^+$-closed.  Thus, we can enumerate the maximal antichains of $\mathbb{P}_{\text{tail}}$, and build a descending sequence of conditions meeting them all through the dense set which has ${\le}\kappa^+$-closure, using this closure and the fact that $M[G_{\kappa}][g]$ is closed under $\kappa^+$-sequences to build through the limit stages.  Then the upward closure of this sequence $G_{\text{tail}} \subseteq \mathbb{P}_{\text{tail}}$ is $M[G_{\kappa}][g]$-generic and is in $V[G_{\kappa}][g]$.  Thus $j(G_{\kappa}) = G_{\kappa} \ast g \ast G_{\text{tail}}$ is $M$-generic, it is in $V[G_{\kappa}][g]$ and $j''G_{\kappa} \subseteq j(G_{\kappa})$.  Thus we have the partial lift $j:V[G_{\kappa}] \to M[j(G_{\kappa})]$.

Next we lift through $\mathbb{Q}$ by finding an $M[j(G_{\kappa})]$-generic filter for $j(\mathbb{Q})$.  The forcing $j(\mathbb{Q})$ adds a club $C \subseteq j(\kappa)$ such that $\forall \delta < j(\kappa)$, $\delta \in C$ implies $o_{\delta^+ \text{-sc}}(\delta)^{M[j(G_{\kappa})]} < j(F)(\delta)$.  Since $j(G_{\kappa}) \in V[G_{\kappa}][g]$ and $M^{\kappa^+} \subseteq M$, it follows that $M[j(G_{\kappa})]^{\kappa^+} \subseteq M[j(G_{\kappa})]$ by [Hamkins1 (Theorem 53)].  Since for every $\beta < \kappa$, there is a dense subset of $\mathbb{Q}$ which is ${\le}\beta$-closed, it follows that there is a dense subset of $j(\mathbb{Q})$ which is ${\le}\kappa^+$-closed.  Also since $|\mathbb{Q}| = \kappa$ and $\mathbb{Q}$ has the $\kappa^+$-chain condition, the forcing $\mathbb{Q}$ has at most $\kappa^+$ many maximal antichains.  Therefore $j(\mathbb{Q})$ has at most $|j(\kappa^+)| \le \kappa^{+^{\kappa^{+^{< \kappa}}}} = \kappa^{++}$ many maximal antichains in $M[j(G_{\kappa})]$.  Thus we can diagonalize to get a generic filter, but we will do so with a master condition as follows.  Let $c = \cup g$ the new club of $\kappa$ and consider $\bar{c} = c \cup \{ \kappa \}$ which is in $M[j(G_{\kappa})]$.  We have $\delta \in c$ implies $o_{\delta^+ \text{-sc}}(\delta)^V < F(\delta)$.  Thus, $\delta \in c$ implies $o_{\delta^+ \text{-sc}}(\delta)^M < F(\delta)$.  By Lemma \ref{Lemma. Mitchell rank for supercompactness doesn't go up.} the Mitchell rank for supercompactness does not go up between $M$ and $M[j(G_{\kappa})]$, so that $\delta \in c$ implies $o_{\delta^+ \text{-sc}}(\delta)^{M[j(G_{\kappa})]} < F(\delta) = j(F)(\delta)$.  Since $o_{\kappa^+ \text{-sc}}(\kappa)^{M[j(G)]} \le o_{\kappa^+ \text{-sc}}(\kappa)^{M} = \beta_0 < \min \{ o_{\kappa^+ \text{-sc}}(\kappa)^V, F(\kappa) \}$ it follows that $o_{\kappa^+ \text{-sc}}(\kappa)^{M[j(G)]} < F(\kappa) = j(F)(\kappa)$.  Thus $\bar{c}$ is a closed, bounded subset of $j(\kappa)$ such that every $\delta$ in $\bar{c}$ has $o_{\delta^+ \text{-sc}}(\delta)^{M[j(G_{\kappa})]} <j(F)(\kappa)$.  Thus $\bar{c} \in j(\mathbb{Q})$.  Thus diagonalize to get an $M[j(G_{\kappa})]$-generic filter $g^* \subseteq j(\mathbb{Q})$ which contains the master condition $\bar{c}$.  Then $j''g \subseteq j(g) = g^*$.  Therefore we may lift the embedding to $j:V[G_{\kappa}][g] \to M[j(G_{\kappa})][j(g)]$.  

The success of the lifting arguments for $o_{\kappa^+ \text{-sc}}(\kappa)^{M} = \beta_0 < \min \{ o_{\kappa^+ \text{-sc}}(\kappa)^V, F(\kappa) \}$, where $\beta_0$ was arbitrary, show that $o_{\kappa^+ \text{-sc}}(\kappa)^{V[G]} \ge \min\{o_{\kappa^+ \text{-sc}}(\kappa)^V, F(\kappa)\}$.  By Lemma \ref{Lemma. Mitchell rank for supercompactness doesn't go up.} we have $o_{\kappa^+ \text{-sc}}(\kappa)^{V[G]} \le o_{\kappa^+ \text{-sc}}(\kappa)^V$.  All that remains is to see that $o_{\kappa^+ \text{-sc}}(\kappa)^{V[G]} \le F(\kappa)$. Consider any $\kappa^+$-supercompact embedding $j:V[G] \to M[j(G)]$ with critical point $\kappa$ in $V[G]$, and the new club $c \subseteq \kappa$, where $\forall \delta < \kappa$, $\delta \in c$ implies $o_{\delta^+ \text{-sc}}(\delta)^V < F(\delta)$.  This club is in any normal fine measure on $P_{\kappa} \kappa^+$.  By the induction hypothesis, $\forall \delta , \kappa$, $\delta \in c$ implies $o_{\delta^+ \text{-sc}}(\delta)^{V[G]} < F(\delta)$.  Applying $j$ to this statement gives $\forall \delta < j(\kappa)$, $\delta \in j(c)$ implies $o_{\delta^+ \text{-sc}}(\delta)^{M[j(G)]} < j(F)(\delta)$.  Since $F \restriction \kappa$ represents $F(\kappa)$, we have $j(F)(\kappa) = F(\kappa)$. Since $\kappa < j(\kappa)$ and $\kappa \in j(c)$ it follows that $o_{\kappa^+ \text{-sc}}(\kappa)^{M[j(G)]} < F(\kappa).$  Thus by Lemma \ref{Lemma. Mitchell rank for supercompactness in M.}, since $j$ was arbitrary, $o_{\kappa^+ \text{-sc}}(\kappa)^{V[G]} \le F(\kappa)$.  Thus $o_{\kappa^+ \text{-sc}}(\kappa)^{V[G]} \le \min \{ o_{\kappa^+ \text{-sc}}(\kappa)^V, F(\kappa) \}$.  Thus for any $\kappa$ which is $\kappa^+$-supercompact for which $F \restriction \kappa$ represents $F(\kappa)$ in V has $o_{\kappa^+ \text{-sc}}(\kappa)^{V[G]} = \min \{ o_{\kappa^+ \text{-sc}}(\kappa)^V, F(\kappa) \}$.  
\end{proof}
\end{theorem}

The most general theorem considers the case of changing Mitchell rank for $\theta$-supercompactness for a $\theta$-supercompact cardinal $\kappa$ when $\theta > \kappa^+$.

\begin{theorem} For any $V \models ZFC + GCH$, any $\Theta : \text{Ord} \to \text{Ord}$ and any $F : \text{Ord} \to \text{Ord}$, there is a forcing extension $V[G]$ where, if $\kappa$ is $\Theta(\kappa)$-supercompact, $\Theta \restriction \kappa$ represents $\Theta(\kappa)$, $F \restriction \kappa$ represents $F(\kappa)$ in V, and $\Theta''\kappa \subseteq \kappa$, then $o_{\Theta(\kappa) \text{-sc}}(\kappa)^{V[g]} =\min\{o_{\Theta(\kappa) \text{-sc}}(\kappa)^V, F(\kappa)\}$.

\begin{proof}
Suppose $V \models ZFC + GCH$. Let $\Theta: \text{Ord} \to \text{Ord}$ and $F : \text{Ord} \to \text{Ord}$.  Let $\mathbb{P}$ be an Easton support Ord-length iteration which forces at $\gamma$ to add a club $c_{\gamma} \subseteq \gamma$ such that $\delta \in c_{\gamma}$ implies $o_{\Theta(\delta) \text{-sc}}(\delta)^V < F(\delta)$, whenever $\gamma$ is a closure point of $\Theta$, meaning $\Theta''\gamma \subseteq \gamma$, otherwise force trivially.  Let $G \subseteq \mathbb{P}$ be $V$-generic.  Let $\delta_0$ be the first inaccessible cardinal.  The forcing before stage $\delta_0$ is trivial, and the forcing at stage $\delta_0$ adds a club to $\delta_0$.  Thus the forcing $\mathbb{P}$ up to and including stage $\delta_0$ has cardinality $\delta_0$ and the forcing past stage $\delta_0$ has a dense set which is closed up to the next inaccessible cardinal.  Thus $\mathbb{P}$ has a closure point at $\delta_0$, thus $V \subseteq V[G]$ satisfies the $\delta_0^+$ approximation and cover properties by [Hamkins4 (Lemma 13)].  Thus by Lemma \ref{Lemma. Mitchell rank for supercompactness doesn't go up.}, the forcing $\mathbb{P}$ does not increase Mitchell rank for $\theta$-supercompactness for cardinals above $\delta_0$ where $\theta$ is also above $\delta_0$.  

Let $\kappa$ be a $\Theta(\kappa)$-supercompact cardinal such that $\Theta''\kappa \subseteq \kappa$, and $\Theta \restriction \kappa$ represents $\Theta(\kappa)$, and  $F \restriction \kappa$ represents $F(\kappa)$ in $V$. Here we will need that representing functions are still representing functions in the extension.  The fact we will use for this is from [Hamkins4] which says that if $V \subseteq V[G]$ satisfies the hypothesis of approximation and cover theorem, then every $\Theta(\kappa)$-closed embedding $j:V[G] \to M[j(G)]$ is a lift of an embedding $j \restriction V : V \to M$ that is in $V$ (and $j \restriction V$ is also a $\Theta(\kappa)$-closed embedding).  Thus, if $F \restriction \kappa$ is a representing function for $F(\kappa)$ in the ground model, $(j \restriction V) (F \restriction \kappa)(\kappa) = F(\kappa)$, and thus $F \restriction \kappa$ is a representing function for $F(\kappa)$ with respect to such embeddings $j$ in $V[G]$.

The induction hypothesis is that any $\kappa' < \kappa$, where $\Theta''\kappa' \subseteq \kappa'$, and $\Theta \restriction \kappa'$ represents $\Theta(\kappa)$, and $F \restriction \kappa'$ represents $F(\kappa')$ in $V$, has $o_{\Theta(\kappa) \text{-sc}}(\kappa')^{V[G]} = \min \{o_{\Theta(\kappa) \text{-sc}}(\kappa')^V, F(\kappa')\}$.  Pick $j:V \to M$ a $\Theta(\kappa)$-supercompactness embedding with critical point $\kappa$ such that $o_{\Theta(\kappa) \text{-sc}}(\kappa)^M = \beta_0 < \min\{ o_{\Theta(\kappa) \text{-sc}}(\kappa)^V, F(\kappa)\}$.  We shall lift $j$ through $\mathbb{P}$.  The forcing above $\kappa$ does not affect the Mitchell rank for $\Theta(\kappa)$-supercompactness of $\kappa$ since the forcing past stage $\kappa$ has a dense subset which is closed up to the next closure point of $\Theta$ past $\kappa$.  Call $\mathbb{P}_{\kappa} \ast \dot{\mathbb{Q}}$ the forcing up to and including stage $\kappa$.  First we shall lift $j$ through $\mathbb{P}_{\kappa}$.  Let $G_{\kappa} \subseteq \mathbb{P}_{\kappa}$ be $V$-generic and $g \subseteq \mathbb{Q}$ be $V[G_{\kappa}]$-generic.  First, we need an $M$-generic filter for $j(\mathbb{P}_{\kappa})$.  Since $cp(j)= \kappa$, the forcings $\mathbb{P}_{\kappa}$ and $j(\mathbb{P}_{\kappa})$ agree up to stage $\kappa$.  The forcing $\mathbb{Q}$ adds a club $c \subseteq \kappa$ such that $\forall \delta < \kappa$, $\delta \in c$ implies $o_{\Theta(\kappa) \text{-sc}}(\delta)^V < F(\delta)$.  Since $o_{\Theta(\kappa) \text{-sc}}(\delta)^V < F(\delta)$ implies $o_{\Theta(\kappa) \text{-sc}}(\delta)^M < F(\delta) = j(F)(\delta)$, and $j(\Theta)(\kappa) = \Theta(\kappa)$, it follows that the $\kappa$th stage of $j(\mathbb{P}_{\kappa})$ is $\mathbb{Q}$.  Thus, $j(\mathbb{P}_{\kappa})$ factors as $\mathbb{P}_{\kappa} \ast \dot{\mathbb{Q}} \ast \mathbb{P}_{\text{tail}}$ where $\mathbb{P}_{\text{tail}}$ is the forcing past stage $\kappa$.  Since $G_{\kappa} \ast g$ is $M$-generic for $\mathbb{P}_{\kappa} \ast \dot{\mathbb{Q}}$, we only need an $M[G_{\kappa}][g]$-generic filter for $\mathbb{P}_{\text{tail}}$ which is $V[G_{\kappa}][g]$.  We will diagonalize to obtain such a filter, after checking that we have met the criteria.  Since $|\mathbb{P}_{\kappa}| = |\mathbb{Q}| = \kappa$, and $V \models \text{GCH}$, it follows that the forcing $\mathbb{P}_{\kappa} \ast \dot{\mathbb{Q}}$ has the $\kappa^+$-chain condition.  Since $M^{\Theta(\kappa)} \subseteq M$, it follows [Hamkins1 (Theorem 54)] that $M[G_{\kappa}][g]^{\Theta(\kappa)} \subseteq M[G_{\kappa}][g]$.  Since $|\mathbb{P}_{\kappa}| = \kappa$, and has the $\kappa^+$-chain condition, it has at most $\kappa^{<\kappa^+} = \kappa^+$ many maximal antichains.  Since $|j(\kappa^+)| \le \kappa^{+^{\Theta(\kappa)^{< \kappa}}} = \Theta(\kappa)^+$, it follows that $\mathbb{P}_{\text{tail}}$ has at most $\Theta(\kappa)^+$ many maximal antichains in $M[G_{\kappa}][g]$.  Since $\mathbb{P}_{\text{tail}}$ has a dense subset which is closed up to the next closure point of $\Theta$ past $\kappa$, it has a dense subset which is ${\le}\Theta(\kappa)$-closed.  Thus, we can enumerate the maximal antichains of $\mathbb{P}_{\text{tail}}$, and build a descending sequence of conditions meeting them all through the dense set which has ${\le}\Theta(\kappa)$-closure, using this closure and the fact that $M[G_{\kappa}][g]$ is closed under $\Theta(\kappa)$-sequences to build through the limit stages.  Then the upward closure of this sequence $G_{\text{tail}} \subseteq \mathbb{P}_{\text{tail}}$ is $M[G_{\kappa}][g]$-generic and is in $V[G_{\kappa}][g]$.  Thus $j(G_{\kappa}) = G_{\kappa} \ast g \ast G_{\text{tail}}$ is $M$-generic, it is in $V[G_{\kappa}][g]$ and $j''G_{\kappa} \subseteq j(G_{\kappa})$.  Thus we have the partial lift $j:V[G_{\kappa}] \to M[j(G_{\kappa})]$.

Next we lift through $\mathbb{Q}$ by finding an $M[j(G_{\kappa})]$-generic filter for $j(\mathbb{Q})$.  The forcing $j(\mathbb{Q})$ adds a club $C \subseteq j(\kappa)$ such that $\forall \delta < j(\kappa)$, $\delta \in C$ implies $o_{\Theta(\kappa) \text{-sc}}(\delta)^{M[j(G_{\kappa})]} < j(F)(\delta)$.  Since $j(G_{\kappa}) \in V[G_{\kappa}][g]$ and $M^{\Theta(\kappa)} \subseteq M$, it follows that $M[j(G_{\kappa})]^{\Theta(\kappa)} \subseteq M[j(G_{\kappa})]$ by [Hamkins1 (Theorem 53)].  Since for every $\beta < \kappa$, there is a dense subset of $\mathbb{Q}$ which is ${\le}\beta$-closed, it follows that there is a dense subset of $j(\mathbb{Q})$ which is ${\le}\Theta(\kappa)$-closed.  Also since $|\mathbb{Q}| = \kappa$ and $\mathbb{Q}$ has the $\kappa^+$-chain condition, the forcing $\mathbb{Q}$ has at most $\kappa^+$ many maximal antichains.  Therefore $j(\mathbb{Q})$ has at most $|j(\kappa^+)| \le \kappa^{+^{\Theta(\kappa)^{< \kappa}}} = \Theta(\kappa)^+$ many maximal antichains in $M[j(G_{\kappa})]$.  Thus we can diagonalize to get a generic filter, but we will do so with a master condition as follows.  Let $c = \cup g$ the new club of $\kappa$ and consider $\bar{c} = c \cup \{ \kappa \}$ which is in $M[j(G_{\kappa})]$.  We have $\delta \in c$ implies $o_{\Theta(\kappa) \text{-sc}}(\delta)^V < F(\delta)$.  Thus, $\delta \in c$ implies $o_{\Theta(\kappa) \text{-sc}}(\delta)^M < F(\delta)$.  By Lemma \ref{Lemma. Mitchell rank for supercompactness doesn't go up.} it follows $\delta \in c$ implies $o_{\Theta(\kappa) \text{-sc}}(\delta)^{M[j(G_{\kappa})]} \le o_{\Theta(\kappa) \text{-sc}}(\delta)^M < F(\delta) = j(F)(\delta)$.  Since $o_{\Theta(\kappa) \text{-sc}}(\kappa)^{M[j(G)]} \le o_{\Theta(\kappa) \text{-sc}}(\kappa)^{M} = \beta_0 < \min \{ o_{\Theta(\kappa) \text{-sc}}(\kappa)^V, F(\kappa) \}$, and $j(\Theta)(\kappa) = \Theta(\kappa)$ it follows that $o_{j(\Theta)(\kappa) \text{-sc}}(\kappa)^{M[j(G)]} < F(\kappa) = j(F)(\kappa)$.  Thus $\bar{c}$ is a closed, bounded subset of $j(\kappa)$ such that every $\delta$ in $\bar{c}$ has $o_{j(\Theta)(\kappa) \text{-sc}}(\delta)^{M[j(G_{\kappa})]} <j(F)(\kappa)$.  Thus $\bar{c} \in j(\mathbb{Q})$.  Thus diagonalize to get an $M[j(G_{\kappa})]$-generic filter $g^* \subseteq j(\mathbb{Q})$ which contains the master condition $\bar{c}$.  Then $j''g \subseteq j(g) = g^*$.  Therefore we may lift the embedding to $j:V[G_{\kappa}][g] \to M[j(G_{\kappa})][j(g)]$.  

The success of the lifting arguments for $o_{\Theta(\kappa) \text{-sc}}(\kappa)^{M} = \beta_0 < \min \{ o_{\Theta(\kappa) \text{-sc}}(\kappa)^V, F(\kappa) \}$, where $\beta_0$ was arbitrary, shows that $o_{\Theta(\kappa) \text{-sc}}(\kappa)^{V[G]} \ge \min\{o_{\Theta(\kappa) \text{-sc}}(\kappa)^V, F(\kappa)\}$.  By Lemma \ref{Lemma. Mitchell rank for supercompactness doesn't go up.} we have $o_{\Theta(\kappa) \text{-sc}}(\kappa)^{V[G]} \le o_{\Theta(\kappa) \text{-sc}}(\kappa)^V$.  All that remains is to see that $o_{\Theta(\kappa) \text{-sc}}(\kappa)^{V[G]} \le F(\kappa)$. Consider any $\Theta(\kappa)$-supercompactness embedding $j:V[G] \to M[j(G)]$ with critical point $\kappa$ in $V[G]$, and the new club $c \subseteq \kappa$, where $\forall \delta < \kappa$, $\delta \in c$ implies $o_{\Theta(\delta) \text{-sc}}(\delta)^V < F(\delta)$.  This club is in any normal fine measure on $P_{\kappa} \Theta(\kappa)$.  By the induction hypothesis, $\forall \delta < \kappa$, $\delta \in c$ implies $o_{\Theta(\delta) \text{-sc}}(\delta)^{V[G]} < F(\delta)$.  Applying $j$ to this statement gives $\forall \delta < j(\kappa)$, $\delta \in j(c)$ implies $o_{j(\Theta)(\kappa) \text{-sc}}(\delta)^{M[j(G)]} < j(F)(\delta)$.  Since $F \restriction \kappa$ represents $F(\kappa)$, we have $j(F)(\kappa) = F(\kappa)$ and since $\Theta \restriction \kappa$ represents $\Theta(\kappa)$, we have $j(\Theta)(\kappa) = \Theta(\kappa)$. Since $\kappa < j(\kappa)$ and $\kappa \in j(c)$ it follows that $o_{\Theta(\kappa) \text{-sc}}(\kappa)^{M[j(G)]} < F(\kappa).$  Thus by Lemma \ref{Lemma. Mitchell rank for supercompactness in M.}, since $j$ was arbitrary, $o_{\Theta(\kappa) \text{-sc}}(\kappa)^{V[G]} \le F(\kappa)$.  Thus $o_{\Theta(\kappa) \text{-sc}}(\kappa)^{V[G]} \le \min \{ o_{\Theta(\kappa) \text{-sc}}(\kappa)^V, F(\kappa) \}$.  Thus for any $\kappa$ which is $\Theta(\kappa)$-supercompact for which $\Theta''\kappa \subseteq \kappa$, and $F \restriction \kappa$ represents $F(\kappa)$, and $\Theta \restriction \kappa$ represents $\Theta(\kappa)$ in V has $o_{\Theta(\kappa) \text{-sc}}(\kappa)^{V[G]} = \min \{ o_{\Theta(\kappa) \text{-sc}}(\kappa)^V, F(\kappa) \}$. 
\end{proof}
\end{theorem}

The last theorem is about large cardinal properties close to the top.  A cardinal $\kappa$ is \textit{huge} with target $\lambda$ if there is an elementary embedding $j:V \to M$, with critical point $\kappa$, such that $j(\kappa) = \lambda$ and $M^{\lambda} \subseteq M$.  A cardinal $\kappa$ is \textit{superhuge} if it is huge with target $\lambda$ for unboundedly many cardinals.

\begin{theorem} If $\kappa$ is huge with target $\lambda$, then there is a forcing extension where $\kappa$ is still huge with target $\lambda$, but $\kappa$ is not superhuge.
\begin{proof} 
Suppose $\kappa$ is huge with target $\lambda$.  Force with Add$(\omega, 1) \ast$ Add$(\lambda^+, 1)$. Then, by Hamkins and Shelah [HS], in the extension, the cardinal $\kappa$ is not $\eta$-supercompact, for all $\eta > \lambda$.  Thus, in the extension, $\kappa$ is not superhuge, since it is not even $\lambda^+$-supercompact, hence not $\lambda^+$-huge.  But, $\kappa$ is still huge with target $\lambda$ since the forcing is $\le \lambda$-closed. 
\end{proof}
\end{theorem}

Open Questions:\\

\noindent 1. How high does the hierarchy of inaccessible degrees go?  Can we define $E_0$-inaccessible?\\
2. If $\kappa$ is a $\Sigma_{n+1}$-reflecting cardinal, is there $V[G]$ where $\kappa$ is $\Sigma_{n}$-reflecting, but not $\Sigma_{n+1}$-reflecting?\\
3.  Is there a forcing like this for $n$-extendible cardinals?\\
4.  Given any large cardinal degree, is there a forcing to softly kill?\\
5. If $o(\kappa)=n$ is there a forcing extension where $\kappa$ is least with $o(\kappa) = n$?\\
6.  Is there a forcing extension where the least strongly compact cardinal is least with $o(\kappa) = 3$, the second strongly
compact cardinal is least above this with $o(\kappa) = 2$ and the third strongly compact cardinal is the least measurable above this one?\\
7.  Is there a forcing extension where the least strongly compact cardinal is the least measurable, the second strongly compact
cardinal is least with $o(\kappa) = 2$, the third strongly compact cardinal is least with $o(\kappa) = 3$, ...\\

\end{document}